\documentclass[aop,MSNbibl,dvips]{arximspdf}

\doi{10.1214/09-AOP507}
\volume{38}
\issue{3}
\pubyear{2010}
\firstpage{1143}
\lastpage{1179}

\makeatletter

\newcommand{\eqref}[1]{(\ref{#1})}

\newtheorem{theorem}{Theorem}
\newtheorem{lemma}{Lemma}
\newtheorem{proposition}{Proposition}
\newproclaim{definition}{Definition}
\newtheorem{corollary}{Corollary}
\newproclaim{remark}{Remark}
\newproclaim{ass}{Assumption}

\makeatother

\begin{document}
\begin{frontmatter}

\title{The obstacle problem for quasilinear stochastic~PDE's\protect\thanksref{TITL}}
\runtitle{Obstacle problem for stochastic PDE's}
\thankstext{TITL}{Supported in part by the research project MATPYL
of the F\'{e}d\'{e}ration de Math\'{e}matiques des Pays de la Loire.}

\begin{aug}
\author[A]{\fnms{Anis} \snm{Matoussi}\corref{}\ead[label=e1]{anis.matoussi@univ-lemans.fr}\thanksref{t1}}
\thankstext{t1}{Supported in part by chaire Risques Financiers de la
fondation du risque,
CMAP-\'{E}cole Polytechniques, Palaiseau-France.}
and
\author[B]{\fnms{Lucretiu} \snm{Stoica}\ead[label=e2]{lstoica@fmi.unibuc.ro}\thanksref{t2}}
\thankstext{t2}{Supported in part by the Contract CEx06-11-18.}
\runauthor{A. Matoussi and L. Stoica}
\affiliation{University of Le Mans and University of Bucharest}
\address[A]{Laboratoire Manceau de Math\'{e}matiques \\
Universit\'{e} du Maine\\ Avenue Olivier Messiaen\\
72085 Le Mans Cedex 9\\
France\\
\printead{e1}} 
\address[B]{Institute of Mathematics ``Simion Stoilow''\\
\quad of the Romanian Academy\\
and \\
Faculty of Mathematics\\
University of Bucharest \\
Str. Academiei 14 \\
Bucharest RO-70109\\
Romania\\
\printead{e2}}
\end{aug}

\received{\smonth{5} \syear{2009}}
\revised{\smonth{9} \syear{2009}}

%
\begin{abstract}
We prove an existence and uniqueness result for the obstacle problem of
quasilinear
parabolic stochastic PDEs. The method is based on the probabilistic
interpretation of the solution
by using the backward doubly stochastic differential equation.
\end{abstract}

%


\end{frontmatter}

\section{Introduction}\label{sec1}
We consider the following stochastic PDE, in $\mathbb{R}^d$,
\begin{eqnarray}
\label{SPDE1}
&&
du_t (x) + \bigl[ \tfrac{1}{2} \Delta u_t (x)  + f_t(x,u_t (x),\nabla u_t (x))\nonumber\\
&&\hspace*{44pt}\hspace*{28pt}{}+ \operatorname{div} g_t (x,u_t ( x )
,\nabla u_t ( x ) ) \bigr] \,dt\\
&&\qquad{} + h_t(x,u_t(x),\nabla u_t(x))\cdot\overleftarrow
{dB}_t = 0,\nonumber
\end{eqnarray}
over the time interval $[0,T]$, with a given final condition $u_T =
\Phi$ and $f, g = (g_1,\ldots,g_d )$,
$h = (h_1,\ldots,h_{d^1} )$ nonlinear random functions. The
differential term with $\overleftarrow{dB}_t$
refers to the backward stochastic integral with respect to a
$d^1$-dimensional Brownian motion on
$ (\Omega, \mathcal{F},\mathbb{P}, (B_t)_{t\geq0} )$. We use the
backward notation because in the proof we will employ the doubly
stochastic framework introduced by Pardoux and Peng \cite
{PardouxPeng94} (see also Bally and Matoussi \cite{BM01} and
Matoussi and Xu \cite{MX2}).

In the case where $f$ and $g$ do not depend of $u$ and $\nabla u$, and
if $h$ is identically null, the equation \eqref{SPDE1} becomes a
linear parabolic equation,
\begin{equation}
\label{LPDE1} \partial_t u (t,x) +
\tfrac{1}{2} \Delta u (t,x) + f(t,x)+ \operatorname{div} g(t,x) =0.
\end{equation}
If $ v \dvtx [0,T] \times\mathbb{R}^d \rightarrow\mathbb{R} $ is a
given function such that $ v (T,x) \leq\Phi(x)$, we may roughly say
that the solution of the obstacle problem for \eqref{LPDE1} is a
function $ u \in\mathbf{L}^2 ([0,T]; H^1 (\mathbb{R}^d) )$ such that
the following conditions are satisfied in $(0,T) \times\mathbb{R}^d$:\vspace*{2pt}
\begin{eqnarray}
\label{OPDE}
\mbox{(i)}&&\qquad\hspace*{-2pt} u \geq v ,\qquad dt\otimes dx  \mbox{-a.e.}, \nonumber\\
\mbox{(ii)}&&\qquad\hspace*{-2pt}   \partial_t u +
\tfrac{1}{2} \Delta u + f + \operatorname{div} g \leq0,
\nonumber\\[-8pt]\\[-8pt]
\mbox{(iii)}&&\qquad\hspace*{-2pt}   (u - v ) \bigl( \partial_t u +
\tfrac{1}{2} \Delta u + f + \operatorname{div} g \bigr) = 0 ,\nonumber\\
\mbox{(iv)}&&\qquad\hspace*{-2pt}   u_T = \Phi, \qquad dx\mbox{-a.e.}\nonumber
\end{eqnarray}
The relation (ii) means that the distribution appearing in the LHS of
the inequality is a nonpositive measure. The relation
(iii) is not rigourously stated. We may roughly say that one has
$\partial_t u +
\frac{1}{2} \Delta u + f + \operatorname{div} g = 0$ on the set $ \{ u > v \}$.

If one expresses the obstacle problem for \eqref{LPDE1} in terms of
variational inequalities one should also ask that the solution has a
minimality property (see Bensoussan--Lions \cite{BensoussanLions78}, page 250, or
Mignot--Puel \cite{MignotPuel75}).

The work of El Karoui et al. \cite{Elk2} treats the obstacle problem
for \eqref{LPDE1} within the framework of backward stochastic
differential equations (BSDE in short). Namely, the equation \eqref
{LPDE1} is considered with $f$ depending of $u$ and $\nabla u$, while
the function $g$ is null (as well $h$) and the obstacle $v$ is
continuous. The solution is represented stochastically as a process and
the main new object of this BSDE framework is a continuous increasing
process that controls the set $\{u=v\}$. This increasing process
determines in fact the measure from the relation (ii). Bally et al.
\cite{BCEF} point out that the continuity of this process allows one
to extend the classical notion of strong variational solution (see
Theorem 2.2 of \cite{BensoussanLions78}, page 238) and express the
solution to the obstacle as a pair $(u, \nu)$ where $ \nu$ equals the
LHS of (ii) and is supported by the set $\{u=v\}$. Moreover, based on
this observation Matoussi and Xu \cite{MX1} generalized the
work under monotonicity and general growth conditions. They have also used
the penalization method and stochastic flow technics (see \cite{BM01}
and \cite{MS03} for
more details on this method). In the present paper, we similarly consider
the solution as a pair $ ( u,\nu),$ point of view which has the
advantage of expressing the notion of solution independently of the
double stochastic framework and without the minimality property of
Mignot--Puel \cite{MignotPuel75}, which would be very difficult to
manipulate in the case of the stochastic PDE.
In Section \ref{measure}, we are going to examine the potential and
the measure associated to a continuous increasing process. We call such
potentials and measures, regular potentials, respectively regular measures.\looseness=1

Now let us consider the final condition to be a fixed function $\Phi
\in\mathbf{L}^2 (\mathbb{R}^d )$ and the obstacle $v$ be a random
continuous function, $ v \dvtx \Omega\times[0,T] \times\mathbb{ R}^d
\rightarrow\mathbb{R }$. Then the obstacle problem for the
equation \eqref{SPDE1} is defined as a pair $ (u, \nu)$, where $ \nu
$ is a random regular measure and $ u \in\mathbf{L}^2 (\Omega\times
[0,T]; H^1 (\mathbb R^d) )$ satisfies the following\ relations:
\begin{eqnarray}\label{OSPDE1}
(\mbox{i}')&&\qquad\hspace*{-4pt} u \geq v, \qquad d\mathbb{P}\otimes dt\otimes dx  \mbox{-a.e.}, \nonumber\\
(\mbox{ii}')&&\qquad\hspace*{-4pt}  du_t (x) + \bigl[
\tfrac{1}{2} \Delta u_t (x) + f_t(x,u_t (x),\nabla u_t (x))\nonumber\\
&&\qquad\hspace*{21pt}\hphantom{du_t (x) + \bigl[}
{}+ \operatorname{div} g_t (x,u_t ( x )
,\nabla u_t ( x ) ) \bigr] \,dt\nonumber\\[-8pt]\\[-8pt]
&&\qquad\hspace*{-4pt}\qquad {} + h_t(x,u_t(x),\nabla u_t(x))\cdot\overleftarrow
{dB}_t = - \nu(dt,dx)\qquad \mbox{a.s.}, \nonumber\\
(\mbox{iii}')&&\qquad\hspace*{-4pt}  \nu(u > v ) =0\qquad \mbox{a.s.},\nonumber\\
(\mbox{iv}')&&\qquad\hspace*{-4pt}  u_T = \Phi,\qquad d\mathbb{P}\otimes dx\mbox{-a.e.}\nonumber
\end{eqnarray}
In Section \ref{section:Ito}, we explain the rigorous sense of the
relation (iii$'$) which is based on the quasi-continuity of $u$. The
main result of our paper is Theorem \ref{maintheorem} which ensures
the existence and uniqueness of the solution of the obstacle problem
for \eqref{SPDE1}. The method of proof is based on the penalization
procedure and the doubly stochastic calculus which is essential,
although the definition of the solution and the statement of the result
avoids the doubly stochastic framework.

Similarly to the case treated in El Karoui et al. \cite{Elk2}, the
most difficult point is to show that the approximating sequence
converges uniformly on the trajectories over the coincidence set $\{u
=v\}$. This is proven in Lemma \ref{essentiel}. The existence and
uniqueness of the solution for equation \eqref{SPDE1} (without
obstacle) has already been proven in \cite{Denis2}. An essential
ingredient in the
treatment of the quasilinear part is the probabilistic representation
of the divergence term obtained in \cite{Stoica} as
well as the doubly stochastic representation corresponding to the
divergence term of the stochastic PDE in \cite{Denis2}.
We must mention the work of Nualard and
Pardoux \cite{Nualart} and Donati-Martin and Pardoux \cite{Donati-Pardoux}
who studied a particular class of obstacle problem  for stochastic PDE driven by
some space--time white noise by using a different techniques.

Finally, we would like to thank our friend Vlad Bally for a stimulating
discussion on the
obstacle problem we had ``la Gare de Montparnasse'' and the referee for
helping us to improve
the presentation.

\section{Preliminaries}
\label{preliminary}\label{spaces}

The basic Hilbert space of our framework is $\mathbf{L}^2 ( \mathbb
{{R}}^d ) $, and we employ the usual notation for its scalar product
and its norm,%
\[
( u,v ) =\int_{\mathbb{R}^d}u ( x ) v (
x )\, dx,\qquad \| u \| _2= \biggl( \int_{\mathbb{R}^d}u^2 (
x ) \,dx \biggr) ^{1/2}.
\]
In general, we shall use the
notation
\[
(u,v)=\int_{\mathbb{R}^d} u(x)v(x)\, dx,
\]
where $u$, $v$ are measurable functions defined in $\mathbb{R }^d$ and
$uv \in\mathbf{L}^1 (\mathbb{R}^d )$.

Our evolution problem will be considered over a fixed time interval
$[0,T]$ and the norm for a function $\mathbf{L}^2 ( [0,T] \times\mathbb{{R}}^d
) $ will be denoted by
\[
\| u \| _{2,2}= \biggl(\int_0^T \int_{\mathbb{R}^d} |u (t,x)|^2 \,dx \,dt
\biggr)^{1/2}.
\]

Another Hilbert space that we use is the first order Sobolev space $H^1
(\mathbb{R }^d%
)= H_0^1 (\mathbb{R }^d%
) .$ Its natural scalar product and norm are%
\[
( u,v ) _{H^1 ( {\mathbb{R}^d} ) }= ( u,v ) + ( \nabla u, \nabla v
),\qquad
\| u \| _{H^1 (
{\mathbb{R}^d} ) }= ( \| u \| _2^2+ \| \nabla
u \| _2^2 ) ^{1/2},
\]
where we denote the gradient by $\nabla u (t,x) = (\partial_1 u (t,x),
\ldots, \partial_d u (t,x) )$.

Of special interest is the subspace $\widetilde{F} \subset\mathbf
{L}^2 ( [0,T]; H^1 ( {\mathbb{R}^d} ) ) $ consisting of all functions
$u(t,x)$ such that $ t \mapsto u_t = u(t, \cdot)$ is continuous in
$\mathbf{L}^2 (\mathbb{R}^d)$. The natural norm on $\widetilde{F}$ is
\[
\| u \| _{T}= \sup_{ 0 \leq t \leq T} \| u_t \|_2 + \biggl(\int_0^T \|
\nabla u_t \|_2 \,dt \biggr)^{1/2}.
\]

The Lebesgue measure in $\mathbb{R}^d$ will be sometimes denoted by
$m$. The space of test functions which we employ in the definition of
weak solutions of the evolution equations \eqref{SPDE1} or \eqref
{LPDE1} is $ \mathcal{D}_T = \mathcal{C}^{\infty} ([0,T] ) \otimes
\mathcal{C}_c^{\infty} (\mathbb{R}^d )$, where
$\mathcal{C}^{\infty} ([0,T] )$ denotes the space of real functions
which can be extended as infinite differentiable functions in the
neighborhood of $[0,T]$ and $ \mathcal{C}_c^{\infty} (\mathbb{R}^d
)$ is the space of infinite differentiable\vspace*{1pt} functions with compact
support in $\mathbb{R}^d$.

\subsection{The probabilistic interpretation of the divergence term}
The operator $ \partial_t + \frac{1}{2} \Delta$, which represents
the main linear part in the equation \eqref{SPDE1}, is
probabilistically interpreted by the Brownian motion in $\mathbb{R}^d$.
We shall view the Brownian motion as a Markov process, and therefore we
next introduce some detailed notation for it. The sample space is $
\Omega' = \mathcal{C } ([0, \infty); \mathbb{R}^d )$, the canonical
process $(W_t)_{t \geq0}$ is defined by $ W_t (\omega) = \omega(t)$,
for any $ \omega\in\Omega'$, $t \geq0$ and the shift operator, $
\theta_t \dvtx \Omega' \rightarrow\Omega'$, is defined by $ \theta
_t (\omega) (s) = \omega(t+s)$, for any $s \geq0$ and $ t \geq0$.
The canonical filtration $ \mathcal{F}_t^0 = \sigma( W_s; s \leq t
)$ is completed by the standard procedure with respect to the
probability measures produced by the transition function
\[
P_t (x, dy) = q_t (x-y)\, dy,\qquad t >0,\ x \in\mathbb{R}^d,
\]
where $ q_t (x) = (2\pi t )^{- d/2} \exp( - |x|^2/2t )$ is
the Gaussian density. Thus, we get a continuous Hunt process $ (\Omega
', W_t, \theta_t, \mathcal{F}, \mathcal{F}_t, \mathbb{P}^x )$. We
shall also use the backward filtration of the future events $ \mathcal
{F}'_t = \sigma(W_s; s \geq t )$ for $t\geq0$. $\mathbb{P}^0$ is
the Wiener measure, which is supported by the set $ \Omega'_0 = \{
\omega\in\Omega', w (0) =0 \}$. We also set $ \Pi_0 (\omega) (t) =
\omega(t) - \omega(0), t \geq0$, which defines a map $ \Pi_0 \dvtx
\Omega' \rightarrow\Omega'_0$. Then $\Pi= (W_0, \Pi_0 ) \dvtx \Omega'
\rightarrow\mathbb{R}^d \times\Omega'_0$ is a bijection. For each
probability measure on $\mathbb{R}^d$, the probability $\mathbb
{P}^{\mu}$ of the Brownian motion started with the initial
distribution $\mu$ is given by
\[
\mathbb{P}^{\mu} = \Pi^{-1} (\mu\otimes\mathbb{P}^0 ).
\]
In particular, for the Lebesgue measure in $\mathbb{R}^d$, which we
denote by $ m = dx$, we have
\[
\mathbb{P}^{m} = \Pi^{-1} (dx\otimes\mathbb{P}^0 ).
\]
These relations are saying that $W_0$ is independent of $ \Pi_0$. It
is known that each component $(W^i_t)_{t \geq0}$ of the Brownian
motion, $ i =1,\ldots,d$, is a martingale under any of the measures
$\mathbb{P}^{\mu}$. The next lemma shows that
$ (W_{t-r}^i, \mathcal{F}'_{t-r} )$, $r \in(0, t]$, is a backward
local martingale under $\mathbb{P}^{m}$.

\begin{lemma}
Let $ 0< s < t$. If $A \in\sigma(W_t)$ is such that $
\mathbb{E}^m [|W_t|; A ] < \infty$, then one has $ \mathbb{E}^m
[|W_s|; A ] < \infty$. Moreover, for each $ B \in\mathcal{F}'_{t}$,
and $ i = 1,\ldots, d$, one has
\[
\mathbb{E}^m [W^i_s; A \cap B ] = \mathbb{E}^m [W^i_t; A \cap B ].
\]
\end{lemma}

\begin{pf}
We note that $W_t$ is uniformly distributed, and
consequently for each $c>0$, the set
$A_c = \{ |W_t| \leq c \}$ satisfies
\[
\mathbb{E}^m [|W_t|; A_c ] < \infty.
\]
This shows that the class of the sets to which applies the statement is
rather large.

The vector
$ (W_0, W_s - W_0, W_t - W_s )$ has the distribution $ m \otimes
\mathcal{N} (0, s) \otimes\mathcal{N} (0, t-s)$, under the measure
$\mathbb{P}^m.$ Then one deduce that $ ( W_s, W_t - W_s )$ has the distribution
$ m \otimes\mathcal{N} (0, t-s)$ and we may write, for $ \varphi_1,
\varphi_2 \in\mathcal{C}_c (\mathbb{R}^d )$,
\begin{eqnarray*}
\mathbb{E}^m [\varphi_1 (W_t -W_s) \varphi_2 (W_t) ] & = &\int
_{\mathbb{R}^d}\int_{\mathbb{R}^d} \varphi_1 (y) \varphi_2 (x+y)
q_{t-s} (y) \,dy \,dx\\
& = &\biggl(\int_{\mathbb{R}^d} \varphi_2 (x) \,dx \biggr) \biggl( \int_{\mathbb{R}^d}
\varphi_1 (y) q_{t-s} (y)\,dy \biggr).
\end{eqnarray*}
This relation shows that the vector $ (W_t - W_s, W_t )$ has the
distribution $\mathcal{N} (0, t-s) \otimes
m $, under $\mathbb{P}^m$.Then the obvious inequality $ |W_s| \leq
|W_t| + |W_t - W_s| (\textbf{1}_{\{|W_t| \leq1\}} + |W_t| )$ allows
one to deduce the first assertion of the lemma.

In order to check the second assertion of the lemma, we write
\[
\mathbb{E}^m [W^i_s; A \cap B ] = \mathbb{E}^m [W^i_t; A \cap B ] -
\mathbb{E}^m [W^i_t - W^i_s; A \cap B ]
\]
and all that it remains to check is that the last term is null. In
order to show this, one first observes that the distribution of the
vector $ (W_t- W_s, W_t , W_{t^1}-W_t, W_{t^2} - W_{t^1},\ldots,
W_{t^n} - W_{t^{n-1}} )$ is $ \mathcal{N} (0, t-s) \otimes m \otimes
\mathcal{N} (0, t^1- t)\otimes\cdots\otimes
\mathcal{N} (0, t^n-t^{n-1})$, for each system $ s <t <t^1<\cdots
<t^n$. Then one has, for each $ B \in\sigma(W_{t^1}-W_t,\ldots,
W_{t^n} - W_{t^{n-1}} )$,
\[
\mathbb{E}^m [W^i_t - W^i_s; A \cap B ] = \mathbb{E}^0 [W^i_t - W^i_s
] m(A) \mathbb{P}^0 (B) = 0,
\]
which implies the assertion of the lemma.
\end{pf}

Now let us assume that $f$ and $|g|$ belong to $ \mathbf{L}^2 ([0,T]
\times\mathbb{R}^d )$ and $ u \in\widetilde{F}$ is a solution of
the deterministic equation \eqref{LPDE1}. Let us denote by
\begin{equation}
\label{representation:divg}
\int_s^t g_r * dW_r = \sum_{i=1}^d \biggl( \int_s^t g_i (r, W_r) \,dW_r^i +
\int_s^t g_i (r, W_r) \,\overleftarrow{dW}{}_r^i \biggr).
\end{equation}
Then one has the following representation (Theorem 3.2 in \cite{Stoica}).

\begin{theorem}
The following relation holds $\mathbb{P}^m$-a.s. for each $
0 \leq s \leq t \leq T$:
\begin{equation}
\label{divergence:term}
\quad u_t (W_t) - u_s(W_s) = \sum_{i=1}^d \int_s^t \partial_i u_r (W_r)
\,dW_r^i - \int_s^t f_r (W_r) \,dr - \frac{1}{2}\int_s^t g_r * dW_r.\hspace*{-12pt}
\end{equation}
\end{theorem}

In \cite{Stoica}, one uses the backward martingale $\overleftarrow
{M}{}^{\mu,i}$ defined under an arbitrary $\mathbb{P}^{\mu}$, with
$\mu$ a probability
measure in $\mathbb{R}^d$, in order to express the integral $\int_s^t
g_r * dW_r$. Though formally the definition looks different, one easily
sees that it is the same object.

\subsection{Regular measures}\label{measure}

In this section, we shall be concerned with some facts related to the
time--space Brownian motion, with the state space $ [ 0,T [ \times
\mathbb{R}%
^{d},$ corresponding to the generator $\partial_{t}+ \frac
{1}{2}\Delta.$ Its associated
semigroup will be denoted by $ ( \widetilde{P}_{t} ) _{t>0}.$ We
may express it in terms of the Gaussian density of the semigroup $ (
P_{t} )
_{t>0}$ in the following way:
\[
\widetilde{P}_{t}\psi( s,x ) = \cases{\displaystyle
\int_{\mathbb{R}^{d}}q_{t} ( x,y ) \psi( s+t,y ) \,dy, &\quad
if $s+t<T$, \vspace*{4pt}\cr
0, &\quad otherwise,%
}
\]
where $\psi\dvtx [ 0,T [ \times\mathbb{R}^{d}\rightarrow\mathbb{R}$
is a bounded Borel measurable function, $s\in[ 0,T [ ,x\in
\mathbb{R}^{d}$ and $t>0.$ So we may also write $ (\widetilde{P}_t
\psi)_s = P_t \psi_{t+s}$ if $ s+t < T$. The corresponding resolvent
has a density
expressed in terms of the density $q_{t}$ too, as follows:
\[
\widetilde{U}_{\alpha}\psi( t,x ) =\int_{t}^{T}\int_{\mathbb{R}%
^{d}}e^{-\alpha( s-t ) }q_{s-t} ( x-y ) \psi(
s,y )\, dy\,ds
\]
or
\[
(\widetilde{U}_{\alpha}\psi)_t =\int_{t}^{T} e^{-\alpha( s-t ) }
P_{s-t} \psi_{s} \,ds.
\]
In particular, this ensures that the excessive functions with respect
to the
time--space Brownian motion are lower semicontinuous. In fact, we will
not use
directly the time space process, but only its semigroup and resolvent. For
related facts concerning excessive functions, the reader is referred to
\cite{BlumenthalGetoor} or \cite{DM}. Some further properties of this semigroup
are presented in the next lemma.

\begin{lemma}
The semigroup $ ( \widetilde{P}_{t} ) _{t>0}$ acts as a strongly
continuous semigroup of contractions on the spaces $\mathbf{L}^{2} ( [
0,T%
[ \times\mathbb{R}^{d} ) =\mathbf{L}^{2} ( [ 0,T [
;\mathbf{L}^{2} ( \mathbb{R}^{d} ) ) $ and $\mathbf{L}^{2} ( [ 0,T%
[ ;H^{1} ( \mathbb{R}^{d} ) ) .$
\end{lemma}

\begin{pf}
Obviously, it is enough to check the following relations:
\begin{eqnarray*}
\lim_{r\rightarrow0} \biggl( \int_{0}^{T-r} \Vert
P_{r}u_{t+r}-u_{t} \Vert_2 ^{2} \,dt+\int_{T-r}^{T} \Vert
u_{t} \Vert_2 ^{2} \,dt \biggr)& =&0,
\\
\lim_{r\rightarrow0} \biggl( \int_{0}^{T-r}\|\nabla(
P_{r}u_{t+r}-u_{t} )\|_2^2 \,dt+\int_{T-r}^{T} \|\nabla u_{t}
\|_2^2 \,dt \biggr) &=& 0.
\end{eqnarray*}
First, we note that for each function $u\in\mathbf{L}^{2} ( [ 0,T [
\times\mathbb{R}^{d} ) $ and $r>0,$ one has%
\[
\lim_{r\rightarrow0}\int_{0}^{T-r} \Vert u_{t+r}-u_{t} \Vert_2
^{2}\,dt=0.
\]
This property is obvious for a function $u\in\mathcal{C}_{c} ( [
0,T [ \times\mathbb{R}^{d} ) $ and then it is obtained by
approximation for any function in $\mathbf{L}^{2} ( [ 0,T [ \times
\mathbb{R}^{d} ) .$ Then the relation%
\[
\lim_{r\rightarrow0}\int_{0}^{T-r} \Vert P_{r}u_{t+r}-u_{t} \Vert_2
^{2}\,dt=0,
\]
easily follows. From it, one deduces the strong continuity of $ (
\widetilde{P}_{t} ) _{t>0}$ on $\mathbf{L}^{2} ( [ 0,T [ \times
\mathbb{R}^{d} ) .$

In order to prove the same property in the space $\mathbf{L}^{2} ( [
0,T%
[ ;H^{1} ( \mathbb{R}^{d} ) ) $, one should start with
the relation%
\[
\lim_{r\rightarrow0}\int_{0}^{T-r} \|\nabla( u_{t+r}-u_{t} ) \|_2^2
\,dt=0,
\]
which holds for each $u\in\mathcal{C}_c^{\infty} ( [ 0,T [ \times
\mathbb{R}^{d} )$ and then repeat, with obvious modifications, the
previous reasoning.
\end{pf}

The next definition restricts our attention to potentials belonging to $
\widetilde{F},$ which is the class of potentials appearing in our parabolic
case of the obstacle problem.

\begin{definition}
\label{potential} (i) A function $\psi\dvtx [ 0,T ] \times\mathbb{R}%
^{d}\rightarrow\overline{\mathbb{R}}$ is called quasicontinuous provided
that for each $\varepsilon>0,$ there exists an open set,
$D_{\varepsilon
}\subset[ 0,T ] \times\mathbb{R}^{d},$ such that $\psi$ is
finite and continuous on $D_{\varepsilon}^{c}$ and%
\[
\mathbb{P}^{m} \bigl( \{ \omega\in\Omega^{\prime}|\exists t\in[ 0,T%
]  \mbox{ s.t. } ( t,W_{t} ( \omega) ) \in D_{\varepsilon
} \} \bigr) <\varepsilon.
\]

(ii) A function $u\dvtx [ 0,T ] \times\mathbb{R}^{d}\rightarrow[
0,\infty] $ is called a regular potential, provided that its
restriction to $[0,T[\times\mathbb{ R}^d$ is
excessive with respect to the time--space semigroup, it is
quasicontinuous, $u\in\widetilde{F}$ and $\lim_{t\rightarrow
T}u_{t}=0$ in $%
\mathbf{L}^{2} ( \mathbb{R}^{d} ) $.
\end{definition}

Observe that if a function $\psi$ is quasicontinuous, then the
process\break
$%
( \psi_{t} ( W_{t} ) ) _{t\in[ 0,T ] }$ is
continuous. Next, we will present the basic properties of the regular
potentials. Do to the expression of the semigroup $ ( \widetilde{P}%
_{t} ) _{t>0}$ in terms of the density, it follows that two excessive
functions which represent the same element in $\widetilde{F}$ should
coincide.

\begin{theorem}
\label{potentielregulier}
Let $u\in\widetilde{F}.$ Then $u$ has a version which is a regular
potential if and only if there exists a continuous increasing process $%
A= ( A_{t} ) _{t\in\lbrack0,T]}$ which is $ ( \mathcal{F}%
_{t} ) _{t\in\lbrack0,T]}$-adapted and such that $%
A_{0}=0$, $\mathbb{E}^{m} [ A_{T}^{2} ] <\infty$ and
\[
\mathrm{(i)}\qquad u_{t}(W_{t})=\mathbb{E} [ A_{T} |\mathcal{F}_{t} ]
-A_{t}\qquad \mathbb{P}^{m}\mbox{-a.s.} 
\]
for each $t\in[ 0,T ] .$ The process $A$ is uniquely determined
by these properties. Moreover, the following relations hold:
\begin{eqnarray*}
\mathrm{(ii)}&&\qquad\hspace*{-4pt} u_{t} ( W_{t} ) =A_{T}-A_{t}-\sum_{i=1}^{d}\int_{t}^{T}\partial
_{i}u_{s} ( W_{s} )\, dW_{s}^{i}\qquad \mathbb{P}^{m}\mbox{-a.s.} , 
\\
\mathrm{(iii)}&&\qquad\hspace*{-4pt}\Vert u_{t} \Vert_2^{2}+\int_{t}^{T} \|\nabla u_{s} \|_2^2
\,ds=\mathbb{E}^{m} ( A_{T}-A_{t} ) ^{2} , 
\\
\mathrm{(iv)}&&\qquad\hspace*{-4pt}( u_{0},\mathcal{\varphi}_{0} ) +\int_{0}^{T} \biggl(\frac{1}{2} ( %
\nabla u_{s}, \nabla\mathcal{\varphi}_{s} ) + ( u_{s},\partial_{s}%
\mathcal{\varphi}_{s} )\biggr)\, ds\\
&&\qquad\hspace*{-4pt}\qquad=\int_{0}^{T}\int_{\mathbb{R}^{d}}%
\mathcal{\varphi} ( s,x ) \nu( ds\,dx )  
\end{eqnarray*}
for each test function $\mathcal{\varphi\in D}_T,$ where $\nu$ is the
measure defined by%
\[
\mathrm{(v)}\qquad\mathcal{\nu} ( \mathcal{\varphi} ) =\mathbb{E}^{m}\int
_{0}^{T}\mathcal{%
\varphi} ( t,W_{t} )\, dA_{t},\qquad \varphi\in \mathcal{C}_{c} ( %
[ 0,T ] \times\mathbb{R}^{d} ). 
\]
\end{theorem}

\begin{pf}
We first remark that the uniqueness of the increasing process in the
representation (i) follows from the uniqueness in the Doob--Meyer
decomposition.

Let us now assume that $\overline{u}$ is a regular potential which is a
version of $u.$ We will use an approximation of $\overline{u}$ constructed
with the resolvent. By the resolvent equation, one has%
\[
\alpha\widetilde{U}_{\alpha}\overline{u}=\alpha\widetilde{U}_{0} (
\overline{u}-\alpha\widetilde{U}_{\alpha}\overline{u} ) .
\]
Let us set $f^{n}=n ( \overline{u}-n\widetilde{U}_{n}\overline{u} )
$ and $u^{n}=n\widetilde{U}_{n}\overline{u}=\widetilde{U}_{0}f^{n}.$
Since $%
\overline{u}$ is excessive, one has $f^{n}\geq0$ and $u^{n},n\in
\mathbb{N}%
^{\ast},$ is an increasing sequence of excessive functions with limit $
\overline{u}.$ In fact $u^{n},n\in\mathbb{N}^{\ast},$ are
potentials and
their trajectories are continuous. On the other hand, the trajectories $
t\rightarrow\overline{u}_{t} ( W_{t} ) $ are continuous on $[0,T[$
by the quasi-continuity of $\overline{u}.$ The process $ ( u_{t} (
W_{t} ) ) _{t\in\lbrack0,T[}$ is a super-martingale, and because $%
\lim_{t\rightarrow T}u_{t}=0$ in $\mathbf{L}^2$, it is a potential
and the trajectories have
null limits at~$T$. Therefore, this approximation also holds uniformly
on the
trajectories, on the closed interval $ [ 0,T ],$%
\[
\lim_{n\rightarrow\infty}\sup_{0\leq t\leq T} \vert u_{t}^{n} (
W_{t} ) -\overline{u}_{t} ( W ) \vert=0\qquad \mathbb{P}^{m}\mbox{-a.s.}
\]
The function $u^{n}$ solves the equation $ ( \partial_{t}+L )
u^{n}+f^{n}=0$ with the condition $u_{T}^{n}=0$ and its backward
representation is%
\[
u_{t}^{n} ( W_{t} ) =\int_{t}^{T}f_{s}^{n} ( W_{s} )
\,ds-\sum_{i=1}^{d}\int_{t}^{T}\partial_{i}u_{s}^{n} ( W_{s} )
\,dW_{s}^{i}.
\]
If we set $A_{t}^{n}=\int_{0}^{t}f_{s}^{n} ( W_{s} )\, ds,$ after
conditioning, this representation gives%
\renewcommand{\theequation}{$\ast$}
\begin{equation}\label{ast}
u_{t}^{n} ( W_{t} )
=A_{T}^{n}-A_{t}^{n}-\sum_{i=1}^{d}\int_{t}^{T}\partial_{i}u_{s}^{n} (
W_{s} ) \,dW_{s}^{i}=\mathbb{E}^{m} [ A_{T}^{n}/\mathcal{F}_{t} ]
-A_{t}^{n}.
\end{equation}
In particular, one deduces%
\[
u_{0}^{n} ( W_{0} ) =\mathbb{E}^{m} [ A_{T}^{n}/\mathcal{F}_{0}%
] =A_{T}^{n}-\sum_{i=1}^{d}\int_{0}^{T}\partial_{i}u_{s}^{n} (
X_{s} ) \,dW_{s}^{i}.
\]
Also from the relation (\ref{ast}), it follows that%
\renewcommand{\theequation}{$\ast\ast$}
\begin{eqnarray}\label{astast}
\mathbb{E}^{m} ( A_{T}^{n}-A_{t}^{n} ) ^{2} &=&\mathbb{E}^{m} \Biggl(
u_{t}^{n} (
W_{t} ) +\sum_{i=1}^{d}\int_{t}^{T}\partial_{i}u_{s}^{n} (
W_{s} ) \,dW_{s}^{i} \Biggr) ^{2} \nonumber\\[-8pt]\\[-8pt]
&=& \Vert u_{t}^{n} \Vert_2 ^{2}+\int_{t}^{T} \|
\nabla u_{s}^{n} \|_2^2 \,ds.\nonumber
\end{eqnarray}
A similar relation holds for differences, in particular one has%
\[
\mathbb{E}^{m} ( A_{T}^{n}-A_{T}^{k} ) ^{2}= \Vert
u_{0}^{n}-u_{0}^{k} \Vert^{2}+2\int_{0}^{T} \|\nabla(
u_{s}^{n}-u_{s}^{k} ) \|_2^2 \,ds.
\]

On the other hand, the preceding lemma ensures that $\lim_{\alpha
\rightarrow\infty}\alpha\widetilde{U}_{\alpha}=I ,$\ in the space $
\mathbf{L}^{2} ( [ 0,T [ ;H^{1} ( \mathbb{R}^{d} ) ) ,$
which implies%
\[
\lim_{n\rightarrow0}\int_{0}^{T} \|\nabla( u_{t}^{n}-\overline{u}%
_{t} ) \|_2^2 \,dt=0.
\]
These last relations imply that there exists a limit $%
\lim_{n}A_{T}^{n}=:A_{T}$ in the sense of $\mathbf{L}^{2} (\mathbb{
P}^{m} ) .$

Let us denote by $M^{n}= ( M_{t}^{n} ) _{t\in[ 0,T ]
},M= ( M_{t} ) _{t\in[ 0,T ] }$ the martingales given by
the conditional expectations $M_{t}^{n}=\mathbb{E}^{m} [
A_{T}^{n}/\mathcal{F}%
_{t} ] ,M_{t}=\mathbb{E}^{m} [ A_{T}/\mathcal{F}_{t} ] .$ Then
one has $\lim_{n \to\infty}M^{n}=M,$ in $\mathbf{L}^{2} ( \mathbb
{P}^{m} ) $, and hence%
\[
\lim_{n \rightarrow\infty}\mathbb{E}^{m}\sup_{0\leq t\leq T} \vert
M_{t}^{n}-M_{t} \vert
^{2}=0.
\]
Then the relation $u_{t}^{n} ( W_{t} ) =M_{t}^{n}-A_{t}^{n}$ shows
that the processes $A^{n},n\in\mathbb{N}^{\ast},$ also converge uniformly
on the trajectories to a continuous process $A= ( A_{t} ) _{t\in%
[ 0,T ] }.$ The inequality%
\[
\sup_{0\leq t\leq T} \vert A_{t}^{n}-A_{t} \vert\leq
A_{T}+ \vert A_{T}^{n}-A_{T} \vert
\]
ensures the conditions to pass to the limit and get%
\[
\lim_{n\rightarrow\infty}\mathbb{E}^{m}\sup_{0\leq t\leq T} \vert
A_{t}^{n}-A_{t} \vert^{2}=0.
\]
Passing to the limit in the relations (\ref{ast}) and (\ref{astast}) one deduces the relations
(i), (ii) and (iii).

In order to check the relation (iv) from the statement, we observe that the
relation is fulfilled by the functions $u^{n},$%
\begin{eqnarray*}
( u_{0}^{n},\mathcal{\varphi}_{0} ) + \int_{0}^{T} \biggl( \frac{1}{2} (
\nabla u_{s}^{n},\nabla\mathcal{\varphi}_{s} ) + (
u_{s}^{n},\mathcal{%
\varphi}_{s} ) \biggr)\, ds &=&\int_{0}^{T}\int_{\mathbb{R}^{d}}\mathcal{%
\varphi} ( s,x ) f^{n} ( s,x ) \,ds\,dx \\
&=&\mathbb{E}^{m}\int_{0}^{T}\mathcal{\varphi} ( s,W_{s} ) \,dA_{s}^{n},
\end{eqnarray*}
where $\mathcal{\varphi}$ is arbitrary in $\mathcal{D}_T.$ In order
to get
the relation (iv), it would suffice to pass to the limit with
$n\rightarrow
\infty$ in this relation. The only term which poses problems is the last
one. The uniform convergence on the trajectories implies that, $\mathbb{P}^{m}$-a.s., the measures $dA_{t}^{n}$ weakly converge to $dA_{t}.$
Therefore, one
has%
\[
\lim_{n\rightarrow\infty}\int_{0}^{T}\mathcal{\varphi}_{t} (
W_{t} )\, dA_{t}^{n}=\int_{0}^{T}\mathcal{\varphi}_{t} (
W_{t} )\, dA_{t}\qquad \mathbb{P}^{m}\mbox{-a.s.}
\]
On the other hand, one has%
\[
\biggl\vert\int_{0}^{T}\mathcal{\varphi}_{t} ( W_{t} )
\,dA_{t}^{n} \biggr\vert\leq\sup_{0\leq t\leq T}\mathcal{\varphi}%
_{t}^{2} ( W_{t} ) +A_{T}^{2}+ \vert A_{T}^{n}-A_{T} \vert
^{2}.
\]
By It\^{o}'s formula and Doob's inequality, one has
\begin{eqnarray*}
\mathbb{E}^m \Bigl(\sup_{0 \leq t \leq T} \mathcal{\varphi}^2 (t, W_t) \Bigr)
& \leq&4 \|\mathcal{\varphi}_0\|^2 +
4 \mathbb{E}^m \biggl( \int_0^T |\partial_t \mathcal{\varphi} (t, W_t) |
\,dt \biggr)^2\\
&&{} +
16 \mathbb{E}^m \int_0^T |\nabla\mathcal{\varphi} |^2(t, W_t) \,dt
\\
& &{} + 2 \mathbb{E}^m \biggl( \int_0^T |\Delta\mathcal
{\varphi} | (t, W_t) \,dt \biggr)^2 \\
& \leq&4 \|\mathcal{\varphi}_0\|^2 + 4 T \int_0^T \|\partial_t
\mathcal{\varphi}_t \|_2^2 \,dt + 16
\int_0^T \|\nabla\mathcal{\varphi}_t \|_2^2 \,dt\\
&&{} + 2T \int_0^T \|
\Delta\mathcal{\varphi}_t \|_2^2 \,dt < \infty.
\end{eqnarray*}
The preceding estimate ensures the possibility of passing to the limit and
deducing that%
\[
\lim_{n}\mathbb{E}^{m}\int_{0}^{T}\mathcal{\varphi} ( s,W_{s} )
\,dA_{s}^{n}=\mathbb{E}^{m}\int_{0}^{T}\mathcal{\varphi} ( s,W_{s} ) \,dA_{s},
\]
and thus we obtain the relation (iv).

Let us now consider the converse. Assume that $u\in\widetilde{F}$ and $A$
is a continuous increasing process adapted to $ ( \mathcal{F}%
_{t} ) _{t\in[ 0,T ] }$ and satisfying the relation (i).
In order to simplify the subsequent notation, it is convenient to
extend our
given function by putting $u_{t}=0$ for $t>T.$ Now, we shall show that
\renewcommand{\theequation}{\arabic{equation}}
\setcounter{equation}{6}
\begin{equation}\label{excessive}
P_{r}(u_{t+r})\leq u_{t},\qquad t\in\lbrack0,T],\ r>0.
\end{equation}
By the Markov property, one gets
\begin{eqnarray*}
P_{r}u_{t+r}(W_{t})& = & \mathbb{E}^{W_{t}} [ u_{t+r}(W_{r}) ] =%
\mathbb{E}^{m} [ u_{t+r}(W_{r+t}) |\mathcal{F}_{t} ] \\
& = & \mathbb{E}^{m} \bigl[ \mathbb{E}^{m} [ A_{T} |\mathcal{F}%
_{t+r} ] -A_{t+r} |\mathcal{F}_{t} \bigr] =\mathbb{E}^{m} %
[ A_{T} |\mathcal{F}_{t} ] -A_{t+r},
\end{eqnarray*}
where the last line comes from the relation (i). This shows that
\[
P_{r}u_{t+r}(W_{t})\leq u_{t}(W_{t})\qquad \mathbb{P}^{m}\mbox{-a.s.}
\]
and as the distribution of $W_{t}$ under $\mathbb{P}^{m}$ is $m$, we
deduce the
inequality \eqref{excessive}. Moreover, this inequality shows by iteration
that if $r\leq r^{\prime},$ then%
\begin{equation}\label{supermedian}
P_{r^{\prime}}u_{t+r^{\prime}}\leq P_{r}u_{t+r}.
\end{equation}
By the properties of the semigroup density and since $t\rightarrow
u_{t}$ is
continuous with values in $\mathbf{L}^{2},$ it follows that, for each
$r>0,$ $%
P_{r}u_{t+r}$, $t\in[ 0,T ] ,$ has a continuous version in $ [
0,T ] \times\mathbb{R}^{d}$ defined by%
\[
\overline{u}^{r} ( t,x ) =\int_{\mathbb{R}^{d}}q_{r} (
x,y ) u_{t+r} ( y ) \,dy.
\]
The inequality (\ref{supermedian}) shows in fact that $\overline
{u}^{r}$ is
supermedian with respect to $ ( \widetilde{P}_{t} ) _{t>0}$ and,
because of continuity, in fact it is excessive. Then $\overline{u}%
=\lim_{r\rightarrow0}\overline{u}^{r}$ is also excessive and since $%
\lim_{r\rightarrow0}P_{r}u_{t+r}=u_{t},$ in $\mathbf{L}^{2},$
clearly $\overline{u}$
is a version of $u.$ The process $ ( \overline{u}_{t} ( W_{t} )
) _{t\in[ 0,T ] }$ is a cdlg supermartingale, and more
precisely a potential. By the relation (i), this process admits a continuous
version. It follows that itself is continuous and, as a consequence,
one has
the following convergence, uniformly on the trajectories:%
\[
\lim_{r\rightarrow0}\sup_{0\leq t\leq T} \vert\overline{u}%
_{t}^{r} ( W_{t} ) -\overline{u}_{t} ( W_t ) \vert
=0\qquad \mathbb{P}^{m}\mbox{-a.s.}
\]
On the other hand, by the representation (i) one has%
\[
\mathbb{E}^{m}\sup_{0\leq t\leq T} \vert\overline{u}_{t} ( W_{t} )
\vert^{2}<\infty,
\]
which leads to%
\[
\lim_{r\rightarrow0}\mathbb{E}^{m}\sup_{0\leq t\leq T} \vert
\overline{u}%
_{t}^{r} ( W_{t} ) -\overline{u}_{t} ( W ) \vert
^{2}=0.
\]
This relation implies that $\overline{u}$ is quasicontinuous, and
hence it
is a regular potential, completing the proof.
\end{pf}

It is known in the probabilistic potential theory that the regular
potentials are associated to continous additive functionals (see \cite
{BlumenthalGetoor},
Section IV.3 or \cite{fot}, Theorem~5.4.2). In the above
theorem, the additive aspect is not evident. In fact, it is hidden in the
relation (i) of Theorem \ref{potentielregulier}. This relation implies that, for $t\leq s,$ $A_{s}-A_{t}$ is
measurable with respect to the completion of $\sigma( W_{r}/r\in[ t,s ]
).$ This can directly be proven but it also follows
from the approximation of $A$ by $A^{n}.$ For the processes $A^{n},n\in
\mathbb{N},$ this measurability property obviously holds. And this
measurability ensures the fact that $A$ corresponds to an additive
functional for the time--space process, which we are not explicitly
using.

The measure $\nu$ from the theorem, expressed in the relation (v), is also
completely determined by the relation (iv), because the test functions are
dense in $\mathcal{C}_{c} ( [ 0,T ] \times\mathbb{R}%
^{d} ) .$ A natural question now is whether one Radon measure on $ [
0,T ] \times\mathbb{R}^{d}$ can be associated via the relation (iv)
from the theorem to two distinct potentials. The answer is that there is
only one such potential and more precisely it can be directly expressed with
the density $q_{t} ( x,y ) $ in terms of the measure, as one can
see from the next lemma.

\begin{lemma}
Let $u$ be a regular potential and $\nu$ a Radon measure on $ [ 0,T%
] \times\mathbb{R}^{d}$ such that relation \textup{(iv)} holds. Then one has%
\[
( \phi,u_{t} ) =\int_{t}^{T}\int_{\mathbb{R}^{d}} \biggl( \int_{%
\mathbb{R}^{d}}\mathcal{\phi} ( x ) q_{s-t} ( x-y )
\,dx \biggr) \nu( ds\,dy )
\]
for each $\phi\in\mathbf{L}^{2} ( \mathbb{R}^{d} ) $ and $t\in[
0,T%
] .$
\end{lemma}

\begin{pf}
We first remark that the relation (iv) is in fact equivalent to the
following more explicit one%
\[
( u_{t},\mathcal{\varphi}_{t} ) +\int_{t}^{T} \biggl( \frac{1}{2}
( \nabla u_{s},\nabla\mathcal{\varphi}_{s} ) + ( u_{s},\partial_{s}%
\mathcal{\varphi}_{s} ) \biggr)\, ds=\int_{t}^{T}\int_{\mathbb{R}^{d}}%
\mathcal{\varphi} ( s,x ) \nu( ds\,dx ) ,
\]
with any $\mathcal{\varphi}\in\mathcal{D}_T$ and $t\in[ 0,T ] .$

Clearly, it is sufficient to prove the lemma for $\mathcal{\phi}\in
\mathcal{%
C}_{c} ( \mathbb{R}^{d} ) $ such that $\mathcal{\phi}\geq0.$ Then
we set $\psi( s,y ) =\int_{\mathbb{R}^{d}}\phi( x )
q_{s-t} ( x-y ) \,dx,$ for $s\in[ t,T ] $ and $y\in
\mathbb{R}^{d}.$ Then $\psi_{s}=P_{s-t}\phi$ and the map
$s\rightarrow
\psi_{s}$ is in $\mathcal{C}^{1} ( ]t,T];L^{2} ( \mathbb{R}%
^{d} ) ) $ and $\partial_{s}\psi= \frac{1}{2} \Delta\psi_{s}.$ Let
$\eta\in
\mathcal{C}_{c} ( \mathbb{R}_{+} ) $ be a decreasing function such
that $\eta=1$ on the interval $ [ 0,1 ] $ and $\eta=0$ for $x\geq
2.$ Set $\eta_{n} ( x ) =\eta( \frac{ \vert x \vert
}{n} ) ,$ so that $ ( \eta_{n} ) _{n\in\mathbb{N}}$ is an
increasing sequence in $\mathcal{C}_{c} ( \mathbb{R}^{d} ) $ with
limit $\mathbf{1}_{\mathbb{R}^{d}}.$ For each fixed $n$, the function
$\eta_{n}\psi$
can be approximated by convolution with smooth functions and then by test
functions from $\mathcal{D}_T,$ and consequently we may write the relation
(iv) in the form%
\begin{eqnarray*}
&&( u_{t},\eta_{n}\psi_{t} ) +\int_{t}^{T} \biggl( \frac{1}{2} (
\nabla u_{s},\nabla(\eta_{n}\psi_{s} ) ) + ( u_{s},\eta
_{n}\,\partial_{s}\psi_{s} ) \biggr)\, ds\\
&&\qquad=\int_{t}^{T}\int_{\mathbb{R}^{d}}\eta_{n} (
x ) \psi( s,x ) \nu( ds\,dx ) .
\end{eqnarray*}
Then it is easy to see that we may pass to the limit with $n\rightarrow
\infty,$ in this relation too. Then we get%
\[
( u_{t},\psi_{t} ) +\int_{t}^{T} \biggl( \frac{1}{2} (
\nabla u_{s},\nabla\psi_{s} ) + ( u_{s}, \partial_{s}\psi
_{s} ) \biggr)\, ds=\int_{t}^{T}\int_{\mathbb{R}^{d}}\psi( s,x ) \nu( ds\,dx ) ,
\]
which becomes the relation asserted by the lemma, on account of the relation
$\partial_{s}\psi=\frac{1}{2} \Delta\psi_{s}.$
\end{pf}

We now introduce the class of measures which intervene in the notion of
solution to the obstacle problem.

\begin{definition}
A nonnegative Radon measure $\nu$ defined in $ [ 0,T ] \times
\mathbb{R}^{d}$ is called regular provided that there exists a regular
potential $u$ such that the relation (iv) from the above theorem is
satisfied.
\end{definition}

As a consequence of the preceding lemma, we see that the regular
measures are
always represented as in the relation (v) of the theorem, with a certain
increasing process. We also note the following properties of a regular
measure, with the notation from the theorem.

\begin{enumerate}
\item A set $B\in\mathcal{B} ( [ 0,T ] \times\mathbb{R}%
^{d} ) $ satisfies the relation $\nu( B ) =0$ if and only if
$\int_{0}^{T}1_{B} ( t,W_{t} )\, dA_{t}=0$ $\mathbb{P}^{m}$-a.s.

\item If a set $B\in\mathcal{B} ( ] 0,T [ \times\mathbb{R}%
^{d} ) $ is polar, in the sense that%
\[
\mathbb{P}^m \bigl( \{ \omega\in\Omega^{\prime}|\exists t\in[ 0,T ]
, ( t,W_{t} ( \omega) ) \in B \} \bigr) =0,
\]
then $\nu( B ) =0.$

\item If $\psi^{1},\psi^{2}\dvtx [ 0,T ] \times\mathbb{R}%
^{d}\rightarrow\overline{\mathbb{R}}$ are Borel measurable and such
that $%
\psi^{1} (t,x)\geq\psi^{2} (t,x), dt\otimes dx$-a.e., and the
processes $ ( \psi_{t}^{i} ( W_{t} ) )
_{t\in[ 0,T ] },i=1,2$, are a.s. continuous, then one has $\nu
( \psi^{1} < \psi^{2} ) =0.$
\end{enumerate}

\subsection{Hypotheses}
\label{Hypotheses}
Let $ B = (B_t)_{t\geq0}$ be a standard
$d^1$-dimensional Brownian motion on a probability
space $ ( \Omega,\mathcal{F}^B, \mathbb{P } )$. So $B_t =
(B_t^1,\ldots, B_t^{d^1} )$ takes values in $\mathbb{R}^{d^1}$. Over
the time interval $[0,T]$ we define the backward filtration $ (\mathcal
{F}_{s,T}^{B} )_{ s\in[0,T]}$ where $\mathcal{F}^B_{s,T}$ is the
completion in $\mathcal{F}^B$ of $\sigma(B_{r}-B_{s};s\leq r\leq T)$.

We denote by $\mathcal{H}_T$ the space of
$H^1(\mathbb{R}^d)$-valued
predictable and $\mathcal{F}^B_{t,T}$-adapted processes $(u_t)_{0\leq
t \leq T}$ such that the trajectories $ t \rightarrow u_t $ are in
$\widetilde{F}$ a.s. and
\[
\| u \| _{T}^2 < \infty.
\]

In the remainder of this paper, we assume that the final condition
$\Phi$ is a given function in $\mathbf{L}^2 (\mathbb{R}^d)$ and the
functions appearing in equation \eqref{SPDE1},
\begin{eqnarray*}
&&f \dvtx [0,T]\times\Omega\times\mathbb{R}^d \times\mathbb{R}\times
\mathbb{R}^{d}
\rightarrow\mathbb{R}, \\
&&g= (g_1,\ldots,g_d) \dvtx [0,T]\times\Omega\times\mathbb{R}^d \times
\mathbb{R}\times\mathbb{R}^{d} \rightarrow\mathbb{R}^{d}, \\
&&h= (h_1,\ldots,h_{d^1}) \dvtx [0,T]\times\Omega\times\mathbb{R}^d
\times\mathbb{R}\times\mathbb{R}^{d}
\rightarrow\mathbb{R}^{d^1},
\end{eqnarray*}
are random functions predictable with respect to the backward
filtration\break  $ (\mathcal{F}_{t,T}^{B} )_{ t\in[0,T]}$. We set
\[
f ( \cdot
,\cdot,\cdot, 0,0):=f^0,\qquad
g( \cdot,\cdot,\cdot,0,0) :=g^0 = (g_1^0,\ldots,g_d^0)
\]
and
\[
h (\cdot,\cdot, \cdot,0,0):=h^0=(h_1^0,\ldots,h_{d^1}^0)
\]
and assume the following hypotheses.

\renewcommand{\theass}{(H)}
\begin{ass}\label{assH}
There exist nonnegative
constants $C, \alpha, \beta$ such that
\begin{longlist}[(iii)]
\item[(i)]
$ |f( t,\omega,x,y,z) -f( t,\omega, x,y^{\prime},z^{\prime}) |
\leq C (|
y-y^{\prime}| + |z-z^{\prime}| ) $.
\item[(ii)] $ (\sum_{j=1}^{d^1}| h_j( t,\omega,x,y,z)
-h_j( t,\omega,x,y^{\prime},z^{\prime})
|^2 )^{1/2}\leq C | y-y^{\prime}|+ \beta|z-z^{\prime}|$.
\item[(iii)]
$ (\sum_{i=1}^{d}| g_i( t,x,y,z) -g_i( t,\omega
,x,y^{{\prime}
},z^{\prime})
|^2 )^{1/2} \leq C | y-y^{\prime}|+ \alpha|
z-z^{\prime}|. $
\item[(iv)] The contraction property (as in \cite{Denis2}):
$ \alpha+\frac{\beta^{2}}{2} < \frac{1}{2} $.
\end{longlist}
\end{ass}

\renewcommand{\theass}{(HD2)}
\begin{ass}\label{assHD2}
\[
\mathbb{E} ( \| f^0 \|_{2,2}^2+ \|
g^0 \| _{2,2}^2+ \| h^0 \| _{2,2}^2 ) <\infty.
\]
\end{ass}

\renewcommand{\theass}{(HO)}
\begin{ass}\label{assHO}
The obstacle $ v (t,\omega, x)$ is a
predictable random function with respect to the backward filtration $
(\mathcal{F}^B_{t,T} )$. We also assume that $ t \mapsto v(t,\omega,
W_t) $ is $\mathbb{P} \otimes\mathbb{P}^m$-a.s. continuous
on $[0,T]$ and satisfies
\[
\label{growth}
v (T, \cdot) \leq\Phi(\cdot).
\]
\end{ass}

We recall that a usual solution (nonreflected one) of the equation
\eqref{SPDE1}
with final condition $ u_T = \Phi$, is a processus $ u \in\mathcal
{H}_T$ such that for each
test function $\varphi\in\mathcal{D}_T$ and any $ \forall t \in
[0,T]$, we have a.s.
\begin{eqnarray}
\label{weak:SPDE}
&& \int_{t}^{T } \biggl[ ( u_{s},\partial_{s}\varphi_s ) +
\frac{1}{2} ( \nabla u_{s}, \nabla\varphi_s ) + ( g_s, \nabla
\varphi_s ) \biggr]\, ds - (\Phi,
\varphi_T ) +
( u_t, \varphi_t ) \nonumber\\[-8pt]\\[-8pt]
&&\qquad = \int_{t}^{T} ( f_s , \varphi_s )\, ds
+ \int_{t}^{T} ( h_s, \varphi_s )\cdot\overleftarrow{dB}_{s}.\nonumber
\end{eqnarray}
By Theorem 8 in \cite{Denis2}, we have existence and uniqueness of the
solution. Moreover, the solution belongs to $\mathcal{H}_T$. We denote
by $ \mathcal{U }(\Phi, f,g,h)$ this solution.

\begin{remark}
Let $L=\sum_{ij}\partial_{i}a^{ij}\,\partial_{j}$ be an elliptic
operator in
divergence form, with the matrix $a= ( a^{ij} ) \dvtx\mathbb{R}^{d}\rightarrow\mathbb{R}^{d}\times\mathbb{R}^{d}$ being symmetric,
measurable and such that%
\[
\lambda\vert\xi\vert^{2}\leq\sum_{ij}a^{ij} ( x )
\xi^{i}\xi^{j}\leq\Lambda\vert\xi\vert^{2}
\]
for any $x,\xi\in\mathbb{R}^{d}.$ If instead of the operator $\frac
{1}{2}%
\Delta$ in our equation (\ref{SPDE1}), we had the operator $L,$ then the contraction
condition (iv) of Assumption  \ref{assH} would be replaced by $\alpha
+\frac{\beta^{2}}{2}<\lambda$ (this ensures the contraction condition as
formulated in \cite{Denis2}). Then the time change $t\rightarrow\frac
{1}{2\Lambda}%
t^{\prime}$ yields a one to one correspondence between the solutions
$u$ of
the equation%
\[
du_{t}+ [ Lu_{t}+f_{t} ( u_{t},\nabla u_{t} ) + \operatorname{div}
g_{t} ( u_{t},\nabla u_{t} ) ] \,dt+h_{t} ( u_{t},\nabla
u_{t} ) \cdot\overleftarrow{dB}_{t}=0,
\]
over $ [ 0,T ] $ and the solutions $\widehat{u}_{t}=u_{1/(
2\Lambda)t}$ satisfying the equation%
\[
d\widehat{u}_{t}+ \bigl[ \tfrac{1}{2}\Delta\widehat{u}_{t}+\widehat{f}%
_{t} ( \widehat{u}_{t},\nabla\widehat{u}_{t} ) +\operatorname{div}\widehat{%
g}_{t} ( \widehat{u}_{t},\nabla\widehat{u}_{t} ) \bigr]\, dt+%
\widehat{h}_{t} ( \widehat{u}_{t},\nabla\widehat{u}_{t} ) \cdot
\overleftarrow{d\widehat{B}}_{t}=0,
\]
over the interval $ [ 0,2\Lambda T ] ,$ with the transformed
coefficients%
\begin{eqnarray*}
\widehat{f} ( t,x,y,z ) &=&\frac{1}{2\Lambda}f \biggl( \frac{1}{%
2\Lambda}t,x,y,z \biggr) ,\widehat{h} ( t,x,y,z ) =\frac{1}{ (
2\Lambda) ^{1/2}}h \biggl( \frac{1}{2\Lambda}t,x,y,z \biggr) ,
\\
\widehat{g}_{i} ( t,x,y,z ) &=&\frac{1}{2\Lambda} \biggl(
g_{i} \biggl( \frac{1}{2\Lambda}t,x,y,z \biggr) +\sum_{j}a^{ij} ( x )
z_{j}-\Lambda z_{i} \biggr) ,\qquad i=1,\ldots,d,
\end{eqnarray*}
and the transformed Brownian motion $\widehat{B}_{t}= ( 2\Lambda)
^{1/2}B_{1/(2\Lambda)t},t\in[ 0,2\Lambda T ] .$
This can be checked just by direct calculations using the above definition
of a solution. Moreover, if one writes $L$ in the form $Lu=\Lambda
\Delta u-\operatorname{div} ( \gamma\nabla u ) ,$ where $\gamma= ( \gamma
^{ij} ) $ is a matrix with the entries $\gamma^{ij} ( x )
=\Lambda\delta^{ij}-a^{ij} ( x ) $, $i,j=1,\ldots,d,$ then one has%
\[
0\leq\gamma=\Lambda I-a\leq( \Lambda-\lambda) I,
\]
in the sense of the order induced by the cone of nonnegative definite
matrices. This implies that one has%
\[
\vert\gamma( x ) \xi\vert\leq( \Lambda
-\lambda) \vert\xi\vert
\]
for any $x,\xi\in\mathbb{R}^{d}.$ Then it easy to deduce that
$\widehat{g}%
_{t} ( x,y,z ) =\frac{1}{2\Lambda} ( g_{1/(2\Lambda)
t} ( x,y,z ) +\gamma( x ) z ) $ fulfils condition
(iii) of Assumption \ref{assH} with a constant\vspace*{1pt} $\widehat{\alpha
}=\frac{1}{%
2\Lambda} ( \alpha+ ( \Lambda-\lambda) ) .$ On the
other hand, one can see that $\widehat{h}$ satisfies condition (ii)
with\vspace*{-3pt} $%
\widehat{\beta}=\frac{1}{ ( 2\Lambda) ^{1/2}}\beta,$ so
that the condition $\alpha+\frac{\beta^{2}}{2}<\lambda,$ ensures $%
\widehat{\alpha}+\frac{\widehat{\beta}^{2}}{2}<\frac{1}{2},$
which is
condition (iv) of our Assumption \ref{assH}. Therefore, we conclude that
our framework covers the case of an equation that involves an elliptic
operator like $L,$ because the properties of the solution $u$ are immediately
obtained from those of the solution $\widehat{u}.$
\end{remark}

\subsection{Quasi-continuity properties}
\label{section:Ito}
In this section, we are going to prove the quasi-continuity of the
solution of the linear
equation, that is, when $f, g, h$ do not depend of $u$ and $ \nabla u$. To
this end, we
first extend the double stochastic It\^{o}'s formula to our framework. We
start by recalling the following result from \cite{Denis2} (stated for
linear SPDE).

\begin{theorem}
\label{FK}
Let $u\in\mathcal{H}_{T}$ be a solution of the equation
\[
du_{t}+ \tfrac{1}{2} \Delta u_t \,dt + (f_t + \operatorname{div} g_{t} )\,dt +
h_t\,\overleftarrow{ dB}_t = 0,
\]
where $f, g, h $ are predictable processes such that
\[
\mathbb{E } \int_0^T [ \|f_t \|_2^2 + \|g_t \|_2^2+ \|h_t \|_2^2 ] \,dt
< \infty\quad\mbox{and}\quad \|\Phi\|_2^2 < \infty.
\]
Then, for any $0\leq s\leq t\leq T$, one has the following stochastic representation,
$ \mathbb{P}\otimes\mathbb{ P}^m$-a.s.,
\begin{eqnarray}\label{Ito:u}
u ( t,W_{t} ) -u ( s,W_{s} )
&=&\sum_{i}\int_{s}^{t}\partial_{i}u ( r,W_{r} )
\,dW_{r}^{i}-\int_{s}^{t}f_r (W_{r} )\, dr \nonumber\\[-8pt]\\[-8pt]
&&{}- \frac{1}{2}%
\int_{s}^{t}g*dW -\int_{s}^{t}h_r ( W_{r} ) \cdot\overleftarrow{dB}_r.\nonumber
\end{eqnarray}
\end{theorem}

We remark that $ \mathcal{F}_T$ and $\mathcal{ F}^B_{0,T}$ are
independent under $ \mathbb{P}\otimes\mathbb{ P}^m$ and therefore in
the above formula the stochastic integrals with respect to $ dW_t$ and
$\overleftarrow{dW}_t$ act independently of $\mathcal{ F}_{0,T}^B$
and similarly the integral with respect to
$\overleftarrow{dB}_t$ acts independently of $ \mathcal{F}_T$.

In particular, the process $ (u_t (W_t) )_{t\in[0,T]} $ admits a
continuous version which we usually denote by $Y= (Y_t )_{t \in[0,T]}$
and we introduce the notation $ Z_t = \nabla u_t (W_t)$. As a
consequence of this theorem, we have the following result.

\begin{corollary}
\label{FK2}
Under the hypothesis of the preceding theorem, one has the following
stochastic representation for $ u^2$, $\mathbb{P}\otimes\mathbb
{P}^{m}$-a.e., for any $0\leq t\leq T$,
\begin{eqnarray}
\label{Ito:2}
u_t^2 (W_{t} ) -\Phi^2 (W_{T} )& = &2 \int_t^T \biggl[ u_s f_s (W_s) - \frac
{1}{2}|\nabla u_s |^2 (W_s)\nonumber\\
&&\hphantom{2 \int_t^T \biggl[}
{} -
\langle\nabla u_s, g_s \rangle(W_s) + \frac{1}{2}|h_s|^2 (W_s) \biggr]\, ds
\nonumber\\[-8pt]\\[-8pt]
&&{}+\int_{t}^{T} ( u_r g_r )(W_r) *dW_r - 2 \sum_{i}\int_{t}^{T} (u_r
\partial_{i}u_r ) (W_{r} )
\,dW_{r}^{i}\nonumber\\
&&{} + 2 \int_{t}^{T} (u_r h_r ) (W_{r} ) \cdot\overleftarrow
{dB}_r.\nonumber
\end{eqnarray}
Moreover, one has the estimate
\begin{eqnarray}
\label{estimationYZ}
&&\mathbb{E} \mathbb{E}^m \Bigl( \sup_{t \leq s \leq T} |Y_s |^2 \Bigr) +
\mathbb{E } \biggl[
\int_t^T \|\nabla u_s \|_2^2 \,ds \biggr]\nonumber\\[-8pt]\\[-8pt]
&&\qquad \leq c \biggl[ \|\phi\|_2^2 + \mathbb{
E} \int_t^T [ \|f_s\|_2^2 + \|g_s\|_2^2 + \|h_s\|_2^2 ] \,ds \biggr]\nonumber
\end{eqnarray}
for each $ t\in[0,T]$.
\end{corollary}

\begin{remark}
With the notation introduced above, one can write the relation \eqref{Ito:2} as
\begin{eqnarray}
\label{Ito:BSDE}
|Y_t|^2 + \int_t^T |Z_r|^2 \,dr & = & |Y_T|^2 + 2 \int_t^T Y_r f_r (W_r)
\,dr - 2 \int_t^T \langle Z_r, g_r (W_r) \rangle\, dr\nonumber\\
&&{} + \int_t^{T} Y_r
g_r (W_r) *dW_r  - 2 \sum_{i}\int_{t}^{T} Y_r Z_{i,r}
\,dW_{r}^{i}\\
&&{} + 2 \int_{t}^{T} Y_r h_r (W_r) \cdot\overleftarrow{dB}_r
+ \int_t^T |h _r |^2 (W_r) \,dr .\nonumber
\end{eqnarray}
\end{remark}

\begin{pf*}{Proof of Corollary \ref{FK2}}
Assume first that $g$ is uniformly bounded and
belongs to $ (\mathcal{H}_T )^d$, so that $E \int_0^T \| \operatorname{div} g_t \|
_2^2 \,dt < \infty$. Then we may represent the solution in the form
\begin{eqnarray*}
u_t (W_{t} ) -u_s (W_{s} )&=&\sum_{i}\int_{s}^{t}\partial_{i}u_r (W_{r} )
\,dW_{r}^{i}-\int_{s}^{t} [ f_r (W_r) + \operatorname{div} g_r (W_r) ] \,dr\\
&&{} -\int_{s}^{t}h_r (W_{r}) \cdot\overleftarrow{dB}_r.
\end{eqnarray*}
By Lemma 1.3 of \cite{PardouxPeng94}, we may write
\begin{eqnarray*}
u_t^2 (W_{t} ) -u_s^2 (W_{s} )& = & - 2 \int_s^t [ u_r (f_r + \operatorname{div} g_r )
(W_r) - |\nabla u_r |^2 (W_r) -|h_r|^2 (W_r) ] \,dr \\
&&{}+ 2 \sum_{i}\int_{s}^{t} (u_r \,\partial_{i}u_r ) (W_{r} )
\,dW_{r}^{i} - 2 \int_{s}^{t} (u_r h_r ) (W_{r} ) \cdot\overleftarrow{dB}_r.
\end{eqnarray*}
On the other hand, by Lemma 3.1 of \cite{Stoica}, one has
\[
- 2 \int_s^t \operatorname{div} (u_r g_r) (W_r) \,dr = \int_s^t u_r g_r (W_r) * dW_r,
\]
so that the preceding relation immediately leads to the relation \eqref
{Ito:2}. Then the standard calculations of BDSDE involving Young's
inequality, BDG inequality and Gronwall's lemma give the estimate
\eqref{estimationYZ}.

Finally, to obtain the result with general $g$ one proceeds by approximation.
\end{pf*}

In the deterministic case, it was proven in \cite{Stoica} that the
solution of a quasilinear equation has a quasicontinuous version. Here,
we shall prove the same property for the solution of an SPDE as is
stated in the next proposition.

\begin{proposition}\label{quasicontinuit:EDPS}
Under the hypothesis of Theorem \ref{FK}, there exists a function $
\overline{u} \dvtx [0,T]\times\Omega\times\mathbb{R}^d \rightarrow
\mathbb{R}$
which is a quasicontinuous version of $u$, in the sense that for each
$\epsilon>0,$ there exits a predictable
random set $D^{\epsilon} \subset
[0,T] \times\Omega\times\mathbb{R}^d $
such that $\mathbb{P}$-a.s. the section $ D_{\omega
}^{\epsilon}$ is open and $\overline{ u} (\cdot, \omega, \cdot) $ is
continuous on its complement $ (D_{\omega}^{\epsilon} )^c$ and
\[
\mathbb{P}\otimes\mathbb{P}^m \bigl( (\omega, \omega') | \exists t \in
[0,T]\mbox{ s.t. }
(t, \omega, W_t (\omega') ) \in
D^{\epsilon} \bigr) \leq\epsilon.
\]
In particular, the process $ (\overline{u}_t (W_t) )_{t \in[0,T]}$
has continuous trajectories,
$\mathbb{P}\otimes\mathbb{P}^m$-a.s.
\end{proposition}

\begin{pf}
Let us choose $ k \in\mathbb{N} $ with $ k > \frac{d}{2}$, so that the
Sobolev space $H^k(\mathbb{R}^d)$ is continuously imbedded in the
space of
H\"{o}lder continuous functions $ \mathcal{C}^{\gamma} (\mathbb{R}^d)$,
with $ \gamma= 1 + [\frac{d}{2}] - \frac{d}{2}$. We first assume that
$ \phi\in H^k(\mathbb{R}^d)$ and $ f$, $g_1,\ldots, g_d$, $h_1,\ldots,h_{d^1}$ belong to $ \mathbf{L}^2 ([0, T] \times\Omega; H^k
(\mathbb{R}^d) )$. By Theorem 8 in \cite{Denis2}, applied with
respect to the Hilbert space $H^k(\mathbb{R}^d)$, one deduces that the
solution $ u = \mathcal{U} (\Phi,f,g,h )$ has the trajectories $ t
\rightarrow u_t (\omega,\cdot)$ continuous in $H^k(\mathbb{R}^d)$
which implies that they are in $ \mathcal{C} [[0,T] \times \mathbb R^d)$.
On the other hand, we have from \eqref{estimationYZ} the following
general estimate
\[
\mathbb{E} \mathbb{E}^m \Bigl( \sup_{ 0 \leq t \leq T} u (t,W_t )^2 \Bigr)
\leq c
\mathbb{E } \biggl[ \|\Phi\|_2^2 + \int_0^T (\|f_t \|_2^2 + \|g_t \|_2^2+
\|h_t \|_2^2 ) \,dt \biggr].
\]
Now, for general $(\Phi,f,g,h)$, one chooses an approximating sequence
of data $(\Phi^n,f^n,g^n,h^n)$ which are $H^k(\mathbb{R}^d)$-valued
and such that
\begin{eqnarray*}
&&\mathbb{E } \biggl( \| \Phi_n - \Phi_{n+1} \|_2 + \int_0^T [\| f_t^n -
f_t^{n+1} \|_2^2 +\| g_t^n - g_t^{n+1} \|_2^2+
\| h_t^n - h_t^{n+1} \|_2^2
] \,dt \biggr)\\
&&\qquad \leq\frac{1}{2^n}.
\end{eqnarray*}
Let $u^n$ be the sequence of $\mathbb{P}$-a.s. continuous
solutions of the equation associated to
$(\Phi^n,f^n,g^n,h^n)$. Then set $E_n^{\epsilon} = \{ | u^n - u^{n+1}
| > \epsilon\} $ and $ D_k^{\epsilon} = \bigcup_{n\geq k}
E_n^{\epsilon}$. Then we have
\begin{eqnarray*}
&&\epsilon^2 \mathbb{P}\otimes\mathbb{P}^m \bigl( (\omega, \omega') |
\exists t \in[0,T] \mbox{ s.t. }
(t, \omega, W_t (\omega') ) \in
E_n^{\epsilon} \bigr)\\
&&\qquad \leq\mathbb{E }\mathbb{E}^m \Bigl[ \sup_{ 0 \leq t
\leq T} \bigl(u_t^n(W_t)- u_t^{n+1}(W_t) \bigr)^2 \Bigr] \leq\frac{c}{2^n}.
\end{eqnarray*}
Further, one takes $\epsilon= \frac{1}{n^2}$ to get
\[
\mathbb{P}\otimes\mathbb{P}^m \bigl( (\omega, \omega') | \exists t \in
[0,T]\mbox{ s.t. }
(t, \omega, W_t (\omega') ) \in
D_k^{\epsilon} \bigr) \leq\sum_{ n=k}^\infty\frac{c n^4}{2^n}.
\]
This shows the statement.
\end{pf}

We also need the quasicontinuity of the solution associated to a random
regular measure, as stated in the next proposition. We first give the
formal definition of this object.

\begin{definition}
\label{random:regular:measure}
We say that $ u \in\mathcal{H}_T$ is a random regular potential
provided that $ u(\cdot,\omega,\cdot)$ has a version which is a
regular potential, $\mathbb{P}(d\omega)$-a.s. The random
variable $ \nu\dvtx \Omega\rightarrow\mathcal{M} ([0,T] \times
\mathbb{R}^d ) $ with values in the set of regular measures on $ [0,T]
\times\mathbb{R}^d$ is called a regular random measure, provided that
there exits a random regular potential $u$ such that the measure $ \nu
(\omega) (dt\,dx)$ is associated to the regular potential $ u (\cdot,
\omega, \cdot)$, $\mathbb{P}(d\omega)$-a.s.
\end{definition}

The relation between a random measure and its associated random regular
potential is described by the following proposition.
\begin{proposition}
\label{quasicontinuit:bis}
Let $u$ be a random regular potential and $\nu$ be the associated
random regular measure. Let $\overline{u}$ be the excessive version of
$u,$ that is, $\overline{u} ( \cdot,\omega,\cdot) $ is a.s. an $ (
\widetilde{P}_{t} ) _{t>0}$-excessive function which coincides with $%
u ( \cdot,\omega,\cdot) $, $dt\,dx$-a.e. Then we have the
following
properties:
\begin{longlist}[(ii)]
\item[(i)] For each $\varepsilon>0,$ there exists a $ ( \mathcal{F}%
_{t,T}^{B} ) _{t\in[ 0,T ] }$-predictable random set $%
D^{\varepsilon}\subset[ 0,T ] \times\Omega\times\mathbb{R}^{d}
$ such that $P$-a.s. the section $D_{\omega}^{\varepsilon}$ is open
and $%
\overline{u} ( \cdot,\omega,\cdot) $ is continuous on its
complement $ ( D_{\omega}^{\varepsilon} ) ^{c}$ and%
\[
\mathbb{P}\otimes\mathbb{P}^{m} \bigl( ( \omega,\omega^{\prime} )
|\exists
t\in[ 0,T ] \mbox{ s.t. }( t,\omega,W_{t} ( \omega)
) \in D_{\omega}^{\varepsilon} \bigr) \leq\varepsilon.
\]
In particular, the process $ (\overline{u}_t (W_t) )_{t \in[0,T]}$
has continuous trajectories,
$\mathbb{P}\otimes\mathbb{P}^m$-a.s.

\item[(ii)] There exists a continuous increasing process $A= ( A_{t} )
_{t\in[ 0,T ] }$ defined on $\Omega\times\Omega^{\prime}$
such that $A_{s}-A_{t}$ is measurable with respect to the $\mathbb
{P}\otimes\mathbb{P}^{m}$-completion of $\mathcal{F}_{t,T}^{B}\vee\sigma( W_{r}/r\in[ t,s%
] ) $, for any $ 0 \leq s \leq t \leq T$, and such that the following
relations are fulfilled
a.s., with any $\mathcal{\varphi}\in\mathcal{D}$ and $t\in[ 0,T ]
$:%
\begin{longlist}[(a)]
\item[(a)] $( u_{t},\mathcal{\varphi}_{t} ) +\int_{t}^{T} ( \frac{1}{2}%
( \nabla u_{s},\nabla\mathcal{\varphi}_{s} ) + ( u_{s},\partial_{s}%
\mathcal{\varphi}_{s} ) ) \,ds=\int_{t}^{T}\int_{\mathbb{R}^{d}}%
\mathcal{\varphi} ( s,x ) \nu( ds\,dx ) $,
\item[(b)] $u_{t}(W_{t})=\mathbb{E} [ A_{T} |\mathcal{F}_{t}\vee\mathcal
{F}%
_{t,T}^{B} ] -A_{t}$,
\item[(c)] $u_{t} ( W_{t} ) =A_{T}-A_{t}-\sum_{i=1}^{d}\int_{t}^{T}\partial
_{i}u_{s} ( W_{s} ) \,dW_{s}^{i}$,
\item[(d)] $\Vert u_{t} \Vert_2^{2}+ \int_{t}^{T} \| \nabla u_{s} \|_2^2
\,ds=\mathbb{E}^{m} ( A_{T}-A_{t} ) ^{2}$,
\item[(e)] $\mathcal{\nu} ( \mathcal{\varphi} ) =\mathbb{E}^{m}\int
_{0}^{T}\mathcal{\varphi} ( t,W_{t} )\, dA_{t}$.
\end{longlist}
\end{longlist}
\end{proposition}

\begin{pf}
The proof of this proposition results from the approximation procedure used
in the proof of Theorem \ref{potentielregulier}.
\begin{longlist}[(ii)]
\item[(i)] Let $ r > 0$. The process $ \overline{u}^r = (\overline{u}^r_t)_{t \in
[0,T]}$, defined by $ \overline{u}^r_t = P_r u _{t+r},$ has the property
that $ (t,x) \rightarrow\overline{u}^r_t$ is jointly continuous
$\mathbb{P}$-a.s. We also have
\[
\lim_{r \to0} \mathbb{E} \mathbb{E}^m \sup_{ 0 \leq t \leq T} |
\overline{u}^r_t (W_t) - \overline{u}_t (W_t) |^2 = 0,
\]
by the arguments used at the end of the proof of Theorem \ref
{potentielregulier}. This one concludes as in the proof of the preceding
proposition.

\item[(ii)] The construction of the increasing process described in Theorem
\ref{potentielregulier} holds globally for a random regular potential
producing on a.e.
trajectory $ \omega\in\Omega$, the increasing process corresponding
to $ u (\cdot,\omega,\cdot)$.\qed
\end{longlist}
\noqed
\end{pf}

We remark that, taking the expectation of the relation (i.i.d.) of this
proposition one gets
\[
\mathbb{E}\mathbb{E}^{m} ( A_{T}^{2} ) =\mathbb{E} \biggl( \Vert
u_{0} \Vert_2 ^{2}+\int_{0}^{T} \|\nabla u_{t} \|_2^2\, dt \biggr) .
\]
%
\section{Existence and uniqueness of the solution of the obstacle
problem}\label{existenceetunicite}
\subsection{The weak solution}
We now precise the definition of the solution of our obstacle problem.
We recall that the data satisfy the hypotheses of Section \ref{Hypotheses}.
\begin{definition}\label{o-spde}
We say that a pair $(u,\nu)$ is a weak solution of the
obstacle problem for the SPDE
\eqref{SPDE1} associated to $(\Phi,f,g,h,v)$, if:
\begin{longlist}[(iii)]
\item[(i)] $ u \in\mathcal{H}_T $ and $u (t,x)\geq v (t,x)$,
$d\mathbb{P}\otimes dt\otimes dx$ a.e. and $u(T,x)=\Phi
(x)$, $d\mathbb{P}\otimes dx$ a.e.,

\item[(ii)] $\nu$ is a random regular
measure on $(0,T) \times\mathbb{R}^d$,

\item[(iii)] for each $\varphi\in\mathcal{D}_T,$ and $ t \in[0,T]$,
\begin{eqnarray}
\label{weak:RSPDE}
&& \int_{t}^{T } \biggl[ ( u_{s},\partial_{s}\varphi_s ) +
\frac{1}{2} ( \nabla u_{s},\nabla\varphi_s ) \biggr]\, ds - (\Phi,
\varphi_T ) +
( u_t, \varphi_t ) \nonumber\\
&&\qquad = \int_{t}^{T} [ (f_s (u_{s},\nabla u_s ), \varphi_s ) -
( g_s (u_{s},\nabla u_s ), \nabla\varphi_s ) ] \,ds \\
&&\quad\qquad{} + \int_{t}^{T} ( h_s ( u_{s},\nabla u_s ) ,\varphi_s ) \cdot
\overleftarrow{dB}_{s}
+ \int_{t}^{T}\int_{\mathbb{R}^d}\varphi_s(x) \nu(ds,dx),\nonumber
\end{eqnarray}

\item[(iv)] if $\overline{u}$ is a quasicontinuous version of $u,$ then one
has%
\[
\int_{0}^{T}\int_{\mathbb{R}^{d}} \bigl( \overline{u}_{s} ( x )
-v_{s} ( x ) \bigr) \nu( ds\,dx ) =0\qquad \mbox{a.s.}
\]
\end{longlist}
\end{definition}

We note that a given solution $u$ can be written as a sum $ u=
u_1+u_2,$ where $u_1$ satisfies a linear equation $ u_1 = \mathcal{U}
(\Phi, f(u, \nabla u), g(u, \nabla u), h(u, \nabla u) )$ with $f, g,
h$ determined by $u$, while $u_2$ is the random regular potential
corresponding to the measure $\nu$. By Propositions \ref
{quasicontinuit:EDPS} and \ref{quasicontinuit:bis}, the conditions
(ii) and (iii)
imply that the process $u$ always admits a quasicontinuous version, so that
the condition (iv) makes sense. We also note that if $\overline{u}$ is
a quasicontinuous version of $u$, then the trajectories of $W$ do not
visit the set $\{ \overline{u} < v \}$, $\mathbb{P}\otimes\mathbb
{P}^m$-a.s.

Here is the main result of our paper.

\begin{theorem}\label{maintheorem}
Assume that the Assumptions \ref{assH}, \ref{assHD2} and \ref{assHO} hold.
Then there exists a unique weak solution of the obstacle problem for
the SPDE \eqref{SPDE1} associated to $(\Phi,f,g,h,v)$.
\end{theorem}

In order to solve the problem, we will use the backward stochastic
differential equation technics. In fact,
we shall follow the main steps of the second proof in \cite{Elk2},
based on the penalization procedure.

The uniqueness assertion of Theorem \ref{maintheorem} results from the
following
comparison result.

\begin{theorem}
\label{comparaison}
Let $\Phi^{\prime},f^{\prime},v^{\prime}$ be similar to $\Phi,f,v$ and
let $ ( u,\nu) $ be the solution of the obstacle problem
corresponding to $ ( \Phi,f,g,h,v ) $ and $ ( u^{\prime},\nu
^{\prime} ) $ the solution corresponding to $ ( \Phi^{\prime
},f^{\prime},g,h,v^{\prime} ) .$ Assume that the following conditions
hold:

\begin{longlist}[(iii)]
\item[(i)] $\Phi\leq\Phi^{\prime}$, $dx\otimes d\mathbb{P}$-a.e.

\item[(ii)] $f ( u,\nabla u ) \leq f^{\prime} (
u,\nabla u ) $, $dt\,dx\otimes\mathbb{P}$-a.e.

\item[(iii)] $v\leq v^{\prime}$, $dt\,dx\otimes\mathbb{P}$-a.e.
\end{longlist}

Then one has $u\leq u^{\prime}$, $dt\,dx\otimes\mathbb{P}$-a.e.
\end{theorem}

\begin{pf}
The proof is identical to that of the similar result of El Karoui et
al. (\cite{Elk2}, Theorem 4.1).

One starts with the following version of It\^{o}'s formula, written with some
quasicontinuous versions $\overline{u},\overline{u}{}^{\prime}$ of the
solutions $u,u^{\prime}$ in the term involving the regular measures
$\nu,\nu^{\prime},$%
\begin{eqnarray*}
&& \mathbb{E} \Vert( u_{t}-u_{t}^{\prime} ) ^{+} \Vert_2
^{2}+\mathbb{E}\int_{t}^{T} \|\nabla(u_{s}-u_{s}^{\prime} )^{+} \|_2^2
\,ds\\
&&\qquad=\mathbb{E} \Vert( \Phi-\Phi^{\prime} )
^{+} \Vert_2 ^{2} +2\mathbb{E}\int_{t}^{T} \bigl( ( u_{s}-u_{s}^{\prime} ) ^{+},f_{s} (
u_{s},\nabla u_{s} ) -f_{s}^{\prime} ( u_{s}^{\prime},\nabla
u_{s}^{\prime} ) \bigr)\, ds\\
&&\quad\qquad{}+2\mathbb{E}\int_{t}^{T}\int_{\mathbb{R}^{d}} (
\overline{u}_{s}-\overline{u}{}^{\prime}_{s} ) ^{+} ( x )
( \nu-\nu^{\prime} ) ( ds\,dx )\\
&&\quad\qquad
{}+2\mathbb{E}\int_{t}^{T} \bigl( \nabla( u_{s}-u_{s}^{\prime} )
^{+},g_{s} ( u_{s},\nabla u_{s} ) -g_{s} (
u_{s}^{\prime},\nabla u_{s}^{\prime} ) \bigr)
\,ds\\
&&\quad\qquad
{}+\mathbb{E}\int_{t}^{T} \Vert h_{s} ( u_{s},\nabla u_{s} )
-h_{s}( u_{s}^{\prime},\nabla u_{s}^{\prime} )
\Vert_2 ^{2}\,ds.
\end{eqnarray*}
We remark that the inclusion $\{ \overline{u}>\overline{u}{}^{\prime}
\} \subset\{ \overline{u}>v \} \cup\{ v> v'\} \cup\{ v' >
\overline{u}'\} $ and the fact that the set $\{ v > v'\} \cup\{v' >
\overline{u}'\}$ is not visited by $W$, imply that $\nu( \overline{u} >
\overline{u}' ) = 0$, a.s. Therefore,
\[
\int_{t}^{T}\int_{\mathbb{R}^{d}} ( \overline{u}_{s}-\overline{u}'_{s} ) ^{+} ( x ) ( \nu-\nu^{\prime} ) (
ds\,dx ) \leq0\qquad \mbox{a.s.}
\]
and then one concludes the proof by Gronwall's lemma.
\end{pf}

\subsection{Approximation by the penalization method}
\label{section:penalization}
For $n\in
\mathbb{N}$, let $ u^n$ be a solution of the following SPDE
\begin{eqnarray}
\label{SPDE:n}
&& du^n_t (x) +
\tfrac{1}{2} \Delta u^n_t (x)\,dt + f(t,x,u^n_t (x),\nabla u^n_t (x))
\,dt\nonumber\\
&&\qquad {}+ n\bigl(u^n_t (x) -v_t (x)\bigr)^{-} \,dt  + \operatorname{div} ( g (t, x,u^n_t ( x )
,\nabla u^n_t ( x ) ) ) \,dt\\
&&\qquad{}
+ h (t,x,u^n_t(x),\nabla u^n_t(x)) \cdot\overleftarrow{dB}_t = 0\nonumber
\end{eqnarray}
with final condition $u^n_T=\Phi$.

Now set $ f_{n}(t,x,y,z)=f(t,x,y,z)+n(y-v_t (x))^{-}$ and $ \nu^n
(dt,dx) := n (u_t^n (x) - v_t (x) )^{-}\, dt\, dx $.
Clearly for each $n \in\mathbb{N}$, $f_n$ is Lipschitz continuous
in $(y,z)$ uniformly in $(t,x)$ with Lipschitz coefficient $C +n$. For
each $ n\in\mathbb{N}$, Theorem 8 in \cite{Denis2} ensures the
existence and uniqueness of a weak solution $ u^n \in\mathcal{H}_T$
of the SPDE \eqref{SPDE:n} associated
with the data $ (\Phi, f_n, g, h)$. We denote by $ Y_t^n = u^n (t, W_t
) $, $ Z_n = \nabla u^n (t, W_t) $ and $S_t = v (t, W_t)$. We shall
also assume that $u^n$ is quasi-continuous, so that $Y^n$ is $\mathbb
{P}\otimes\mathbb{P}^m$-a.e. continuous. Then $
(Y^n, Z^n )$ solves the BSDE associated to the data $ (\Phi, f_n, g, h)$
\begin{eqnarray}
\label{BSDE:n}
Y_{t}^n& = &\Phi( W_T ) +\int_{t}^{T}f_r
(W_{r},Y^n_{r},Z^n_{r}
) \,dr + n \int_t^T (Y_r^n - S_r^n)^{-} \,dr\nonumber\\
&&{} +\frac{1}{2}\int_{t}^{T}g_r
(W_{r},Y^n_{r},
Z^n_{r} ) *dW\nonumber\\[-8pt]\\[-8pt]
&&{} + \int_{t}^{T}h_r
(W_{r},Y^n_{r},Z^n_{r}
)\cdot\overleftarrow{dB}_r\nonumber\\
&&{} -\sum_{i}\int
_{t}^{T}Z^n_{i,r}\,dW_{r}^{i} .\nonumber
\end{eqnarray}
We define $K_{t}^{n}=n
\int_{0}^{t}(Y_{s}^{n}-S_{s})^{-}\,ds$ and establish the following lemmas.

\begin{lemma}\label{penalization:estimate1}
The triple $(Y^{n},Z^{n}, K^n)$ satisfies the following estimates
\begin{eqnarray}\label{estimate1:n}
&&\mathbb{E} \mathbb{E}^m |Y_t^n |^2 +
\lambda_{\varepsilon} \mathbb{E } \mathbb{E}^m\int_t^T |Z_r^n|^2
\,dr\nonumber\\
&&\qquad\leq c \mathbb{ E}\mathbb{ E}^m \biggl[ |\Phi(W_T) |^2 + \int_t^T \bigl(
|f^0_s (W_s) |^2 + |g^0_s (W_s)|^2 + |h^0_s (W_s)|^2 \bigr)\, ds \biggr]\nonumber\\[-8pt]\\[-8pt]
&&\quad\qquad{} + c_{\varepsilon} \mathbb{E} \mathbb{E}^m \int_t^T |Y_r^n |^2 \,dr +
c_{\delta} \mathbb{E} \mathbb{E}^m
\Bigl(\sup_{t \leq r \leq T} |S_r|^2 \Bigr)\nonumber\\
&&\quad \qquad {} + \delta\mathbb{E} \mathbb{E}^m
(K_T^n - K_t^n )^2,\nonumber
\end{eqnarray}
where $\lambda_{\varepsilon} = 1 - 2 \alpha- \beta^2 - \varepsilon$,
$c_{\varepsilon}, c_{\delta}$ are a positive constants and $\varepsilon>0,
\delta>0 $ can be chosen small enough such that $\lambda_{\varepsilon}
> 0 $.
\end{lemma}

\begin{pf}
By using It\^{o}'s formula \eqref{Ito:BSDE} for $(Y^n,
Z^n),$ we get
\begin{eqnarray}\label{Ito:BSDEn}
|Y_t^n|^2  + \int_t^T |Z_r^n|^2\, dr &=& | \Phi(W_T)|^2 + 2 \int_t^T
Y_s^n f_s ( W_s, Y_s^n, Z_s^n ) \,ds  \nonumber\\
&&{}+ 2 \int_t^T Y_s^n \,dK_s^n - 2 \int_t^T \langle Z_s^n, g_s ( W_s,Y_s^n,Z_s^n ) \rangle \,ds\nonumber\\
&&{} +
\int_t^{T} Y_s^n g_s(W_s, Y_s^n, Z_s^n )*dW - 2 \sum_{i}\int_{t}^{T}
Y_s ^n Z_{i,s}^n\, dW_{s}^{i} \\
&&{} + 2 \int_{t}^{T} Y_s^n h_s(W_s, Y_s^n, Z_s^n ) \cdot\overleftarrow
{dB}_s\nonumber\\
&&{} + \int_t^T |h_s(W_s,Y_s^n, Z_s^n|^2 \,ds .\nonumber
\end{eqnarray}
Using Assumption \ref{assH} and taking the expectation in the above
equation under $ \mathbb{P} \otimes\mathbb{P}^m$, we get
\begin{eqnarray*}
&&\mathbb{E} \mathbb{E}^m | Y_{t}^{n} | ^{2}+ \mathbb{E} \mathbb{E}^m
\int_{t}^{T} | Z_{s}^{n} | ^{2}\,ds \\
&&\qquad \leq \mathbb{E} | \Phi(W_T) |
^{2}+ c_{\varepsilon} \mathbb{E} \mathbb{E}^m \int_{t}^{T} [
|f_s^0(W_s)|^{2}+ |g_s^0(W_s)|^{2} + |h_s^0(W_s)|^{2} ] \,ds\\
&&\quad\qquad{} + c_{\varepsilon} \mathbb{ E} \mathbb{E}^m \int_{t}^{T} |
Y_{s}^{n} |^{2}\,ds + ( 2\alpha+ \beta^2 + \varepsilon) \mathbb{E} \mathbb{E}^m
\int_{t}^{T} | Z_{s}^{n} | ^{2}\,ds\\
&&\quad\qquad{} +\frac{1}{%
\gamma}\mathbb{E} \mathbb{E}^m \Bigl[\sup_{t\leq s \leq
T}|S_{s}|^{2}\Bigr]+\gamma
\mathbb{E}\mathbb{E}^m [(K_{T}^{n}-K_{t}^{n})^{2}],
\end{eqnarray*}
where $\varepsilon>0$, $\gamma>0$ are a arbitrary constants and $
c_{\varepsilon}$ is a constant which can be \mbox{different} from line to line.
We have used the inequality $ \int_t^T Y_s^n \,dK_s^n \geq\int_t^T
S_s^n \,dK_s^n$ and then we have applied Schwartz's inequality. We also
have used the fact that \mbox{under} the measure $\mathbb{P}^{m}$ the
forward--backward integral $\int Y_r^n g (r,W_r, Y_r^n, Z_r^n )*dW $ as
well the other stochastic integrals with respect to the brownian terms
have null expectation under $ \mathbb{P} \otimes\mathbb{P}^m$.
Finally, Gronwall's lemma leads to the desired inequali\-ty.\looseness=1 
\end{pf}

\begin{lemma}
\begin{eqnarray}
\label{estimation:Kn}
\mathbb{E}\mathbb{E}^m [(K_{T}^{n}-K_{t}^{n})^{2}] &\leq& c' [ \mathbb
{E} \mathbb{E}^m  | Y_{t}^{n} | ^{2} + \|\Phi\|_2^2 ]\nonumber\\
&&{} +
c_{\varepsilon} \biggl[ \mathbb{ E} \mathbb{E}^m \int_{t}^{T} [ |
Y_{s}^{n} | ^{2}
+ | Z_{s}^{n} | ^{2} ]\, ds\\
&&\hspace*{25pt}{}  + \mathbb{E} \int_{t}^{T} [ \|f_s^0\|_2^{2}+ \|g_s^0\|_2^{2} + \|
h_s^0\|_2^{2} ]\, ds \biggr].\nonumber
\end{eqnarray}
\end{lemma}

\begin{pf}
Let now $(\widetilde{u}^n)_{n\in\mathbb{N}}$ be the weak solutions
of the following linear type equations
\[
d \widetilde{u}_t^{n} + \tfrac{1}{2} \Delta\widetilde{u}_t^{n} +
\operatorname{div} g_t (u_t^n, \nabla u _t^n ) \,dt + h_t (u_t^n, \nabla u_t^n ) \cdot
\overleftarrow{dB}_t= 0,
\]
with final condition $\widetilde{u}_T^n = 0.$ Set $ \widetilde{
Y}_t^n = \widetilde{u}^n (t, W_t ) $ and $ \widetilde{Z}_n = \nabla
\widetilde{u}^n (t, W_t) $.
Then by the estimate \eqref{estimationYZ}, one has
%
\begin{equation}
\label{estimation:tilde}
\mathbb{E} \mathbb{E}^m \biggl[ | \widetilde{Y}_{t}^{n} | ^{2} + \int_0^T
| \widetilde{Z}_{s}^{n} | \,ds \biggr] \leq\tilde{c} \Lambda,
\end{equation}
where $ \Lambda= \mathbb{E} \mathbb{E}^m \int_{0}^{T} [ |g_s(W_s,
Y_s^n, Z_s^n)|^{2} + |h_s(W_s, Y_s^n, Z_s^n)|^{2} ] \,ds$.
Since $ u^n - \widetilde{u}^n$ verifies the equation
\[
\partial_t(u_t^n - \widetilde{u}_t^{n}) + \tfrac{1}{2} \Delta(u^n -
\widetilde{u}_t^{n}) + f_t(u_t^n, \nabla u_t^n) + n(u^n_t -v_t)^{-}\, dt = 0,
\]
we have the stochastic representation
\begin{eqnarray*}
Y_{t}^n - \widetilde{Y}^n_t & = &\Phi( W_T ) +\int_{t}^{T}f_r
(W_{r},Y^n_{r},Z^n_{r}
) \,dr + K_T^n - K_t^n \\
&&{}-\sum_{i}\int_{t}^{T} (Z^n_{i,r}- \widetilde
{Z}^n_{i,r} )\, dW_{r}^{i}
\end{eqnarray*}
%
from which one easily obtains the estimate
\begin{eqnarray*}
&&\mathbb{E}\mathbb{E}^m [(K_{T}^{n}-K_{t}^{n})^{2}]\\
&&\qquad\leq c \mathbb{E}\mathbb{E}^m \biggl[  | Y_{t}^{n} | ^{2} + |\widetilde{ Y}_{t}^{n}| ^{2} +
|\Phi( W_T ) |^2 \\
&&\quad\qquad\hphantom{c \mathbb{E}\mathbb{E}^m \biggl[}
{}+ \int_t^T ( |f_s^0 (W_s) |^2 + | Y_{s}^{n} | ^{2}+
| Z_{s}^{n} |^{2} ) \,ds
 + \int_t^T |\widetilde{Z}_{s}^{n}|^{2} \,ds \biggr].
\end{eqnarray*}
Hence, using \eqref{estimation:tilde}, we get
\begin{eqnarray*}
&&\mathbb{E}\mathbb{E}^m [(K_{T}^{n}-K_{t}^{n})^{2}] \\
&&\qquad \leq c' \mathbb
{E} \mathbb{E}^m [ | Y_{t}^{n} | ^{2} +
|\Phi(W_T) |^2 ]\\
&&\quad\qquad
{} + c_{\varepsilon}' \mathbb{E} \mathbb{E}^m \biggl[ \int
_t^T ( | Y_{s}^{n} | ^{2} + | Z_{s}^{n} |^{2} ) \,ds \\
&&\quad\qquad\hphantom{{} + c_{\varepsilon}' \mathbb{E} \mathbb{E}^m \biggl[}
{} + \int_t^T [ |f_s^0 (W_s) |^2 + |g_s^0 (W_s) |^2
+ |h_s^0 (W_s) |^2 ]\, ds \biggr],
\end{eqnarray*}
which gives our assertion.
\end{pf}

\begin{lemma}\label{mainestimate}
The triple $(Y^{n},Z^{n},K^{n})$ satisfies the following estimate
\begin{eqnarray*}
&&\mathbb{E} \mathbb{E}^m \Bigl( \sup_{0 \leq s \leq T} |Y_s^n |^2 \Bigr) +
\mathbb{E }\mathbb{E}^m \int_0^T |Z_s^n|^2 \,ds + \mathbb{E} \mathbb
{E}^m ( K_T^n )^2 \\
&&\qquad \leq c \biggl[ \|\Phi\|_2^2 + \mathbb{E} \mathbb{E}^m
\Bigl(\sup_{0 \leq s \leq T} |S_s|^2 \Bigr)
 + \mathbb{E} \int_0^T [ \|f^0_s\|_2^2 + \|g^0_s\|_2^2 + \|h^0_s\|
_2^2 ] \,ds \biggr],
\end{eqnarray*}
where $ c >0$ is a constant.
\end{lemma}

\begin{pf}
From \eqref{estimate1:n} and \eqref{estimation:Kn}, we get
\begin{eqnarray*}
&&(1 - \delta c' ) \mathbb{E} \mathbb{E}^m |Y_s^n |^2 +
 ( 1 - 2 \alpha- \beta^2 - \varepsilon- \delta c_{\varepsilon}' )
\mathbb{E } \mathbb{E}^m\int_s^T |Z_r^n|^2 \,dr \\
&&\qquad\leq( 1 + c' \delta
) \|\Phi\|_2^2 + (c_{\varepsilon} + \delta c_{\varepsilon}' ) \Lambda+
 (c_{\varepsilon} + \delta c_{\varepsilon}' ) \mathbb{E} \mathbb{E}^m
\int_s^T |Y_r^n |^2 \,ds\\
&&\quad\qquad{} + c_{\delta} \mathbb{E} \mathbb{E}^m
\Bigl(\sup_{t \leq r \leq T} |S_r|^2 \Bigr),
\end{eqnarray*}
where $ \Lambda= \mathbb{E} \mathbb{E}^m \int_t^T [|f^0_s (W_s) |^2
+ |g^0_s (W_s)|^2 + |h^0_s(W_s)|^2 ] \,ds $. It then follows from
Gronwall's lemma that
\begin{eqnarray*}
&&\sup_{0 \leq s \leq T} \mathbb{E} \mathbb{E}^m ( |Y_s^n |^2 ) +
\mathbb{E }\mathbb{E}^m \int_s^T |Z_r^n|^2 \,dr + \mathbb{E} \mathbb
{E}^m ( K_T^n )^2 \\
&&\qquad \leq c_1 \biggl[ \|\Phi\|_2^2 + \mathbb{E} \mathbb{E}^m
\Bigl(\sup_{0 \leq r \leq T} |S_r|^2 \Bigr)
 + \mathbb{ E}\int_s^T [ \|f^0_r\|_2^2 + \|g^0_r\|_2^2 + \|h^0_r\|
_2^2 ] \,dr \biggr].
\end{eqnarray*}
Coming back to the equation \eqref{BSDE:n} and using
Bukholder--Davis--Gundy inequality and the last estimates, we get our statement.
\end{pf}

In order to prove the strong convergence of the sequence $
(Y^n,Z^n,K^n)$, we shall need the following result.

\begin{lemma}[(The essential step)]
\label{essentiel}
\begin{equation}
\lim_{n \to\infty} \mathbb{E} \mathbb{E}^m \Bigl[
\sup_{ 0\leq t \leq T} \bigl( (Y_t^n - S_t )^- \bigr)^2 \Bigr] = 0.
\end{equation}
\end{lemma}

\begin{pf}
Let $(u^n)_{ n\in\mathbb{N}}$ be the sequence of
solutions of the penalized SPDE defined in \eqref{SPDE:n}. From Lemma
\ref{mainestimate}, it follows that the sequence $ ( f (u^n, \nabla
u^n), g (u^n,\break \nabla u^n), h (u^n, \nabla u^n) )_{n\in\mathbb{N}}$ is
bounded in $\mathbf{L}^2 ( [0,T] \times\Omega\times\mathbb{R}^d ;
\mathbb{R}^{1+d +d^1} )$. We may
choose then a subsequence which is weakly convergent to a system of
predictable processes $(\bar{f},\bar{g},\bar{h})$ and,
on account of the Lemma \ref{Mazur} in the \hyperref[app]{Appendix}, we obtain a
sequence of families of coefficients of convex combinations,
$(a^k)_{k\in\mathbb{N}},$ such that the sequences
\[
\hat{f}^k = \sum_{i \in I_k} \alpha_i^k f (u^i,\nabla u^i ),\qquad
\hat{g}^k = \sum_{i \in I_k} \alpha_i^k g (u^i,\nabla u^i )
\]
and
\[
\hat{h}^k = \sum_{i \in I_k} \alpha_i^k h (u^i,\nabla u^i )
\]
converge strongly, that is,
\[
\lim_{k \to\infty} \mathbb{E }\int_0^T \| \hat{f}_t^k - \bar{f}_t \|_2^2 \,dt = 0
\]
and similarly for $\hat{g}^k$, $\bar{g}$ and $\hat{h}^k$, $\bar{h}$.

Now for $i \geq n$, we denote by $u^{i,n}$ the solution of the equation
\begin{eqnarray}\label{SPDE:i}
&&du_t^{i,n} + \bigl[\tfrac{1}2 \Delta u_t^{i,n} - n u_t^{i,n} + n v_t + f_t
(u^i,\nabla u^i ) + \operatorname{div} g_t (u^i,\nabla u^i ) \bigr] \, dt\nonumber\\[-8pt]\\[-8pt]
&&\qquad{} + h_t (u^i,\nabla u^i ) \cdot\overleftarrow{dB}_t = 0\nonumber
\end{eqnarray}
with final condition $u_T ^{i,n} = v_T$. By comparison (Theorem \ref
{comparaison}), we have that $ u^{i,n} \leq u^i $. Further, we set
$ \hat{u}^k =\sum_{i \in I_k} \alpha_i^k u^{i,n_k},$ where $ n_k=
\inf I_k$ and we deduce that
\begin{equation}
\label{monotone}
\hat{u}^k \leq\sum_{i \in I_k} \alpha_i^k u^{i} \leq\lim_{n \to
\infty}u^n,
\end{equation}
where the last inequality comes from the monotonicity of the sequence $u^n$.
Moreover, we observe that $ \hat{u}^k$ is a solution of the equation
\begin{equation}
\label{SPDE:k}
d\hat{u}_t^k + \bigl[\tfrac{1}{2} \Delta\hat{u}_t^k - n_k \hat{u}_t^k +
n_k v_t + \hat{f}_t^k + \operatorname{div} \hat{g}_t^k \bigr]\, dt +
\hat{h}_t^k \cdot\overleftarrow{dB}_t = 0
\end{equation}
with final condition $\hat{u}_T^k = v_T$.

Now we are going to take the advantage of the fact that the equations
satisfied by the sequence of solutions $\hat{u}^k$ have strongly
convergent coefficients. Let us denote by $\widehat{Y}^k$ the
continuous version on $[0,T]$ of the process
$ (\hat{u}^k (W_t) )_{t \in[0,T]}$, for any $ k \in\mathbb{N}$.
We will prove now that there exists a subsequence such that
\begin{equation}
\label{Ychapeau}
\lim_{ k \to\infty} \sup_{0 \leq t \leq T} |\widehat{Y}^k_t - S_t|
=0\qquad \mathbb{P}\otimes\mathbb{P}^m\mbox{-a.s.}
\end{equation}
Since the equation
\eqref{SPDE:k} is linear, the solution decomposes as a sum of four
terms each corresponding to one of the coefficients $\hat{f}^k, \hat
{g}^k, \hat{ h}^k, v.$ So it is enough to treat separately each term.

(a) In the case where $f \equiv0$, $g \equiv0$, $h \equiv0$ one
obtains the term corresponding to $v$. Then the relation
\eqref{Ychapeau} is a direct consequence of the Lemma \ref{obstacleS}.

(b) In the case where $v \equiv0$, $g \equiv0$, $h \equiv0$, the
representation of $\widehat{Y}^k$ is given by
\[
\label{Ito:Ykf}
\widehat{Y}_t^k = \int_t^T e^{-n_k (s-t)} \hat{f}^k_s (W_s) \,ds -
\sum_{i =1 }^d \int_{t}^{T} e^{-n_k (s-t)}
\,\partial_{i}\hat{u}_r^k (W_{s} )\,dW_{s}^{i}.
\]
Thus, we have
\[
\biggl| \int_t^T e^{-n_k (s-t)} \hat{f}_s^k (W_s) \,ds \biggr| \leq\frac
{1}{\sqrt{2n_k}} \biggl(\int_t^T
(\hat{f}_s^k (W_s) )^2 \,ds \biggr)^{1/2}.
\]
This shows that $ \lim_{k\to\infty} \sup_{0\leq t \leq T} |\int
_t^T e^{-n_k (s-t)} \hat{f}^k_s (W_s)\, ds | =0, \mathbb{P}\otimes
\mathbb{P}^m$-a.s., on some subsequence. For the second term
in the expression of
$\widehat{Y}^k$, we make an integration by parts formula to get
\[
\int_t^T e^{-n_k(s-t)} \,\partial_{i}\hat{u}_s^k (W_{s} )
\,dW_{s}^{i} = e^{-n_k (T-t)}U_T^{i,k} - U_t^{i,k} + n_k \int_t^T
U_s^{i,k} e^{-n_k(s-t)} \,ds,
\]
where $ U_s^{i,k} = \int_0^s \partial_{i} \hat{u}_r^k (W_{r} )
\,dW_{r}^{i}$. By the Corollary \ref{coefficientfn} of Section \ref
{technical:results}, we know that the\vspace*{1pt} martingales $U^{i,k}, k \in
\mathbb{N}$, converges to zero in $\mathbf{L}^2,$ and hence on a
subsequence we have $ \lim_{k\to\infty} \sup_{0 \leq t \leq T}
|U_t^{i,k} | = 0, \mathbb{P}\otimes\mathbb{P}^m$-a.s. Then
by Lemma \ref{mc}, we see that for that subsequence
\[
\lim_{ k \to\infty} \sup_{0\leq t \leq T} \biggl|\int_t^T e^{-n_k(s-t)}
\,\partial_{i}\hat{u}_s^k (W_{s} )
\,dW_{s}^{i} \biggr|=0\qquad \mathbb{P}\otimes\mathbb{P}^m\mbox{-a.s.}
\]
Therefore, the desired result \eqref{Ychapeau} holds also in this
case. This time we get $ \lim_{k\to\infty} \sup_{0 \leq t \leq T}
|\widehat{Y}_t^{k} | = 0, \mathbb{P}\otimes\mathbb{P}^m$-a.s.

(c) In the case where $f \equiv0$, $h \equiv0$, $v \equiv0$, the
representation of $\widehat{Y}^k$ is given by
\begin{eqnarray*}
\widehat{Y}_t^k & = & \int_t^T e^{-n_k (s-t)} \hat{{g}}^k_s * dW -
\sum_{i =1 }^d \int_{t}^{T} e^{-n_k (s-t)}
\,\partial_{i}\hat{u}_r^k (W_{s} )\,dW_{s}^{i}\\
& = & \sum_i \int_t^T e^{-n_k (s-t)} \hat{{g}}^k_s (W_s) \,dW_s^i +
\sum_i \int_t^T e^{-n_k (s-t)} \hat{{g}}^k_s (W_s)\,\overleftarrow{
dW}{}^i_s\\
&&{}- \sum_{i =1 }^d \int_{t}^{T} e^{-n_k (s-t)}
\,\partial_{i}\hat{u}_r^k (W_{s} )\,dW_{s}^{i}.
\end{eqnarray*}
Now the proof is similar to that of the preceding case. We treat only
the second term in the last expression. We set
$ \overleftarrow{U}{}^{i,k}_s = \int_s^T \hat{{g}}^k_r
(W_r)\,\overleftarrow{ dW}{}^{m,i}_r$. Integration by parts formula gives
\[
\int_t^T e^{-n_k(s-t)} \,d\overleftarrow{U}{}^{i,k}_s = \overleftarrow
{U}{}^{i,k}_t - e^{-n_k (T-t)} \overleftarrow{U}{}^{i,k}_T
- n_k \int_t^T \overleftarrow{U}{}^{i,k}_s e^{-n_k(s-t)} \,ds.
\]
On the other hand, the convergence $ \hat{{g}}^k \rightarrow\bar{g}$
implies that the backward martingale $ (\overleftarrow{U}{}^{i,k}_t )_{t
\in[0,T]}$ converges
to $ (\int_t^T \bar{g}_{i,r} (W_r) \,\overleftarrow{ dW}{}^{m,i}_r )_{t
\in[0,T]}$ in $\mathbf{L}^2 (\mathbb{P}\otimes\mathbb{P}^m )$.
The other terms in the above expression of $\widehat{Y}^k$ may be handled
similarly by integration by parts and taking into account Corollary
\ref{coefficientgn}.
Using again Lemma \ref{mc}, as in the preceding case, we get the
relation \eqref{Ychapeau} in the form $ \lim_{k\to\infty} \sup_{0
\leq t \leq T} |\widehat{Y}_t^{k} | = 0, \mathbb{P}\otimes\mathbb
{P}^m$-a.s.

(d) In the case where $f \equiv0$, $g \equiv0$, $v \equiv0$, the
representation of $\widehat{Y}^k$ is given by
\[
\widehat{Y}_t^k  = - \sum_{i =1 }^d \int_{t}^{T} e^{-n_k (s-t)}
\,\partial_{i}\hat{u}_r^k (W_{s} )\,dW_{s}^{i} + \int_t^T e^{-n_k (s-t)}
\hat{{h}}^k_s \cdot\overleftarrow{dB}_s .
\]
On account of Lemma \ref{coefficienthn}, the same arguments used in
the previous cases work again.

Now it is easy to see that the relation \eqref{Ychapeau} holds for the
general case. On the other hand, \eqref{monotone} and \eqref{Ychapeau}
clearly imply the relation
\[
\lim_{n \to\infty}
\sup_{ 0\leq t \leq T} (Y_t^n - S_t )^- = 0 \qquad\mathbb{ P}\otimes
\mathbb{P}^m \mbox{-a.s.}
\]
and then, since $Y^n$ is bounded in $\mathbf{L}^2$, one gets the
relation of our statement.
\end{pf}

We have also the following result.

\begin{lemma}
\label{convergence:YZK}
There exists a progressively measurable triple of processes
$ (Y_t, Z_t, K_t )_{t \in[0,T]}$ such that
\begin{eqnarray}
\label{convergence:YZKn}
&&\mathbb{E} \mathbb{E}^m \biggl[ \sup_{0 \leq s \leq T} |Y_t^n - Y_t|^2
+
\int_0^T |Z_t^n - Z_t|^2 \,dt\nonumber\\[-8pt]\\[-8pt]
&&\hspace*{104pt}{} + \sup_{ 0 \leq t \leq T} |K_t^n - K_t|^2
\biggr] \longrightarrow0\qquad
\mbox{as } n \to\infty.\nonumber
\end{eqnarray}
Moreover we have that $ (Y_t, Z_t, K_t )_{t \in[0,T]}$ satisfies
$Y_t\geq S_t, \forall t \in[0,T] $ and $ \int_{0}^{T}(Y_s-S_s)\,dK_s =
0$, $ \mathbb{P}\otimes\mathbb{P}^m$-a.e.
\end{lemma}

\begin{pf}
From the monotonicity of the sequence $ (f_n)_{n \in
\mathbb{N}}$ and the comparison Theorem \ref{comparaison}, we get
that $ u^n (t,x) \leq u^{n+1} (t,x)$, $ dt\,dx\otimes\mathbb{P}$-a.e., therefore one has $ Y_t^n \leq Y_t^{n+1}$, for all $ t \in
[0,T]$, $\mathbb{P}\otimes\mathbb{P}^m$-a.s. Thus, there exists a
predictable real valued process
$ Y = (Y_t )_{t \in[0,T]}$ such that $ Y_t^n \uparrow Y_t,$ for all $t
\in[0,T]$ a.s. and by Lemma \ref{mainestimate} and Fatou's lemma, one
gets
\[
\mathbb{E} \mathbb{E}^m \Bigl( \sup_{0 \leq s \leq T} |Y_t|^2 \Bigr) \leq c.
\]
Moreover, from the dominated convergence theorem one has
\begin{equation}
\label{convergence:Y}
\mathbb{E} \mathbb{E}^m \int_0^T |Y_t^n - Y_t|^2\, dt \longrightarrow
0\qquad \mbox{as } n \to\infty.
\end{equation}
The relation \eqref{Ito:BSDE} gives, for $n \geq p$,
\begin{eqnarray}
\label{Ito:BSDEn}
&& |Y_t^n -Y_t^p|^2 + \int_t^T |Z_s^n-Z_s^p|_a^2 \,ds\nonumber\\
&&\qquad = 2 \int_t^T
(Y_s^n - Y_s^p ) [f_s ( W_s, Y_s^n, Z_s^n ) - f_s ( W_s, Y_s^p, Z_s^p )
] \,ds \nonumber\\
&&\quad\qquad{} + 2 \int_t^T (Y_s^n - Y_s^p ) \,d (K_s^n - K_s^p )\\
&&\quad\qquad{} - 2 \int_t^T \langle Z_s^n - Z_s^p, g_s ( W_s,Y_s^n,Z_s^n ) -
g_s ( W_s,Y_s^p,Z_s^p ) \rangle \,ds \nonumber\\
&&\quad\qquad{} + \int_t^{T} (Y_s^n - Y_s^p ) [g_s ( X_s, Y_s^n, Z_s^n ) - g_s (
W_s, Y_s^p, Z_s^p ) ] *dW\nonumber \\
&&\quad\qquad{} - 2 \sum_{i}\int_{t}^{T} (Y_s^n - Y_s^p ) (Z_{i,s}^n - Z_{i,s}^p )
\,dW_{s}^{i}\nonumber\\
&&\quad\qquad{}  + 2 \int_{t}^{T} (Y_s^n - Y_s^p ) [ h_s(W_s, Y_s^n, Z_s^n )
- h_s(W_s, Y_s^p, Z_s^p) ] \cdot\overleftarrow{dB}_s \nonumber\\
&&\quad\qquad{} + \int_t^T |h_s(W_s,Y_s^n, Z_s^n )- h_s(W_s,Y_s^p, Z_s^p)|^2 \,ds .\nonumber
\end{eqnarray}
By standard calculation, one deduces that
\begin{eqnarray}\label{cauchy:Z}
\mathbb{E} \mathbb{E}^m \int_t^T |Z_s^n-Z_s^p|^2 \,ds & \leq &c
\mathbb{E} \mathbb{E}^m \int_t^T |Y_s^n - Y_s^p|^2\nonumber\\
&&{}
+ 4 \mathbb{E} \mathbb{E}^m \int_t^T (Y_s^n - S_s )^- \,dK_s^p
\\
&&{} + 4 \mathbb{E} \mathbb{E}^m \int_t^T (Y_s^p - S_s )^- \,dK_s^n.\nonumber
\end{eqnarray}
Therefore from Lemma \ref{essentiel}, \eqref{convergence:Y} and
\eqref{cauchy:Z} one gets
\begin{eqnarray}
\label{cauchy:YZ}
&&\mathbb{E} \mathbb{E}^m \int_0^T |Y_t^n - Y_t^p|^2 \,dt\nonumber\\[-8pt]\\[-8pt]
&&\qquad {} + \mathbb{E}
\mathbb{E}^m \int_0^T |Z_t^n - Z_t^p|^2 \,dt \longrightarrow0
\qquad \mbox{as } n,p \to\infty.\nonumber
\end{eqnarray}
The rest of the proof is the same as in El Karoui et al. (\cite{Elk2},
pages 721--722), in particular we get that there exists a pair $(Z,K)$ of
progressively measurable processes with values in $ \mathbb{R}^d
\times\mathbb{ R}$ such that
\begin{eqnarray}
&&\mathbb{E} \mathbb{E}^m \biggl[ \sup_{0 \leq s \leq T} |Y_t^n - Y_t|^2
+ \int_0^T |Z_t^n - Z_t|^2 \,dt\nonumber\\
&&\hspace*{102pt}{} + \sup_{ 0 \leq t \leq T} |K_t^n - K_t|^2
\biggr] \longrightarrow0
\qquad \mbox{as } n \to\infty.\nonumber
\end{eqnarray}

It is obvious that $ (K_t )_{t \in[0,T]}$ is an increasing continuous process.
On the other hand, since from Lemma \ref{essentiel} we have $\lim
_{n\rightarrow
\infty} \mathbb{E} \mathbb{E}^m [\sup_{0 \leq t \leq
T}((Y^n_t-S_t)^-)^2 ]=0$, then, $\mathbb{P}\otimes\mathbb{P}^m$-a.s.,
\begin{equation}\label{condition:obstacle}
Y_t\geq S_t\qquad \forall t \in[0,T],
\end{equation}
which yields that $\int_{0}^{T}(Y_s-S_s)\,dK_s \geq0$. Finally, we also have
$\int_{0}^{T}(Y_s-S_s)\,dK_s = 0$ since on the other hand the
sequences $(Y^n)_{n\geq0}$ and $(K^n)_{n\geq0}$ converge
uniformly (at least for a subsequence), respectively, to $Y$ and $K$
and
\[
\int_{0}^{T}( Y_{s}^{n}-S_{s})\,dK_{s}^{n}=-n\int_0^{T}\bigl(( Y_{s}^{n}-S_{s}) ^{-}\bigr) ^{2}\,ds\leq0.
\]
\upqed
\end{pf}

As a consequence of the last proof, we obtain the following
generalization of the RBSDE introduced in \cite{Elk2}.
\begin{corollary}
\label{RBDSDE:definition}
The limiting triple of processes $ (Y_t, Z_t, K_t )_{t \in[0,T]}$ is a
solution of the following reflected backward doubly stochastic
differential equation (in short RBDSDE):
\begin{eqnarray}
\label{RBDSDE}
Y_{t}& = &\Phi( W_T ) +\int_{t}^{T}f_r(W_{r},Y_{r},Z_{r}) \,dr + K_T - K_t \nonumber\\
&&{}+\frac{1}{2}\int_{t}^{T}g_r (W_{r},Y_{r},Z_{r} )
*dW\\
&&{} + \int_{t}^{T}h_r(W_{r},Y_{r},Z^n_{r})\cdot\overleftarrow{dB}_r -\sum_{i}\int_{t}^{T}Z_{i,r}\,dW_{r}^{i}\nonumber
\end{eqnarray}
with $Y_t\geq S_t, \forall t \in[0,T] $, $ (K_t )_{t \in[0,T]}$ is
an increasing continuous process, $K_0=0$ and
\begin{equation}
\label{reflection:minimum}
\int_{0}^{T}(Y_s-S_s)\,dK_s = 0.
\end{equation}
\end{corollary}

\begin{pf*}{Proof of Theorem \ref{maintheorem}}
Since
\[
\int_0^T ( \| u_t^n - u_t^p \|_2^2 + \|\nabla u_t^n - \nabla u_t^p \|
_2^2 )\,dt = \mathbb{E}^m \int_0^T ( | Y_t^n - Y_t^p |^2 + | Z_t^n -
Z_t^p |^2 ) \,dt,
\]
by the preceding lemma one deduces that the sequence $ (u^n )_{n\in
\mathbb{N}}$ is a Cauchy sequence in $
\mathbf{L}^2 ( \Omega\times[0,T]; H^1 (\mathbb{R}^d) )$ and hence
has a limit $u$ in this space. Also from the preceding lemma, it
follows that $dK_t^n$ weakly converges to $dK_t$, $\mathbb{P}\otimes
\mathbb{P}^m$-a.e. This implies that
\begin{eqnarray*}
\lim_{n} \int_0^T \int_{\mathbb{R}^d} n ( u^n - v )^- \varphi(t,x)
\,dt\, dx &=& \lim_{n} \mathbb{E}^m \int_0^T \varphi_t(W_t) \,dK_t^n \\
&=& \int_0^T \int_{\mathbb{R}^d} \varphi(t,x) \nu(dt\, dx ),
\end{eqnarray*}
where $\nu$ is the regular measure defined by
\[
\int_0^T \int_{\mathbb{R}^d} \varphi(t,x) \nu(dt\, dx ) = E^m \int
_0^T \varphi_t(W_t)\, dK_t .
\]
Writing the equation \eqref{SPDE:n} in the weak form and passing to
the limit one obtains the equation \eqref{weak:RSPDE} with $u$ and
this $\nu$. The arguments we have explained after Definition \ref
{o-spde} ensure that $u$ admits a quasicontinuous version $\overline{u}$.
Then one deduces that $ ( \overline{u}_t(W_t) )_{t \in[0,T]}$ should
coincide with $(Y_t)_{t \in[0,T]}$, $\mathbb{P}\otimes\mathbb
{P}^m$-a.e. Therefore, the inequality $ Y_t \geq S_t$ implies
$ u \geq v $, $dt \otimes\mathbb{P} \otimes dx$-a.e. and
the relation $ \int_0^T (Y_t - S_t )\, dK_t = 0$ implies the relation
(iv) of Definition \ref{o-spde}.
\end{pf*}

\section{Some technical lemmas}
\label{technical:results}
\begin{lemma}
\label{coefficientf}
Let $f \in\mathbf{L}^2 ( [0,T] \times\mathbb{R}^d ; \mathbb{R} )$
and denote by $(u^n)_{n\in\mathbb{N}}$ the sequence of solutions of
the equations
\[
\bigl(\partial_t + \tfrac{1}{2} \Delta\bigr) u^{n} - n u^{n} + f = 0\qquad \forall
n \in\mathbb{N},
\]
with final condition $u_T^n = 0$. Then we have
\begin{equation}
\label{deterministe:n}
\int_0^T \|\nabla u_t^n \|_2^2 \,dt \leq c \biggl[ \frac{1}{n} \int_0^T \|
f_t\|_2^2 \,dt + \int_0^T e^{-2n (T-t)}
\|f_t\|_2^2 \,dt \biggr].
\end{equation}
\end{lemma}

\begin{pf}
It is well known that the solution $(u^n)_{n\in
\mathbb{ N}}$ is expressed in terms of the semigroup $P_t$ by
\[
u_t^n = \int_t^T e^{-n(s-t)} P_{s-t} f_s \,ds .
\]
A direct calculation shows that one has
\[
n \int_t^T e^{-n(s-t)} P_{s-t}u_s^0 \,ds = u_t^0 - u_t^n,
\]
which leads to
\begin{equation}
\label{u_n}
u_t^n = e^{- n(T-t)}u_t^0 + n \int_t^T e^{-n(s-t)} (u_t^0 - P_{s-t}
u_s^0 ) \,ds.
\end{equation}
The function $ \overline{u}^n_t = e^{- n(T-t)}u_t^0 $ is a solution of the
equation\vspace*{1pt} $
(\partial_t + \frac{1}{2} \Delta) \overline{u}^{n} - n \overline{u}^{n} +
\bar{f} = 0$
where $ \bar{f}_t = e^{(T-t)} f_t$. Therefore, one has the following
estimate for the gradient of the first term
in the expression of $u^n$
\begin{equation}
\label{gradient1}
\int_0^T e^{-n(T-t)} \|\nabla u_t^0 \|_2^2 \,dt \leq c \int_0^T e^{-2n
(T-t)} \|f_t\|_2^2 \,dt
\end{equation}
(see Lemma 5 of \cite{Denis2} for details).
In order to estimate the gradient of the second term of the expression
of $u^n$, we first remark that
\[
u_t^0 - P_{s-t} u_s^0 = \int_t^s P_{r-t} f_r \,dr,
\]
so that one has
\begin{eqnarray*}
\biggl\|n\nabla\int_t^T e^{-n(s-t)} (u_t^0 - P_{s-t} u_s^0 )\, ds \biggr\|_2 &
\leq &  n \int_0^{T-t} e^{-ns}
\int_0^s \|\nabla P_{r} f_{t+r} \|_2 \,dr \,ds \\
& \leq & nc \int_0^{T-t} e^{-ns} \int_0^s \frac{1}{\sqrt{r}}\|
f_{t+r} \|_2 \,dr \,ds,
\end{eqnarray*}
where we have used the well-known inequality
\[
\|\nabla P_{r} \varphi\|_2 \leq\frac{c}{\sqrt{r}} \|\varphi\|_2\qquad
\mbox{for }  \varphi\in\mathbf{L}^2 .
\]
Then we estimate the time integral of the norm of the gradient, which
is the expression we are interested in,
\begin{eqnarray*}
&&\int_0^T \biggl\|n\nabla\int_t^T e^{-n(s-t)} (u_t^0 - P_{s-t} u_s^0 ) \,ds
\biggr\|_2^2 \,dt \\
&&\qquad\leq
c^2 \int_0^T \biggl[\int_0^{T-t} n e^{-ns} \int_0^s \frac{1}{\sqrt{r}}\|
f_{t+r} \|_2 \,dr \,ds \biggr]^2 \,dt\\
&&\qquad =  c^2 \int_0^T \int_0^s \int_0^T \int_0^{s'} \int_0^{T - s \vee
s'} n e^{-ns} n e^{-ns'}\frac{1}{\sqrt{r}}\|f_{t+r} \|_2\\
&&\quad\qquad\hspace*{118pt}{}\times
\frac{1}{\sqrt{r'}}\|f_{t+r'} \|_2 \,dt \,dr' \,ds' \,dr \,ds\\
&&\qquad \leq \int_0^T \|f_t\|_2^2 \,dt \biggl(\int_0^T \frac{1}{2} \sqrt{s} n
e^{-ns}\, ds \biggr)^2 \leq\frac{c}{n} \int_0^T \|f_t\|_2^2 \,dt.
\end{eqnarray*}
This estimate together with \eqref{gradient1} imply the statement
\eqref{deterministe:n}.
\end{pf}

Obviously, the lemma implies that $ \lim_{n\to\infty} \int_0^T \|
\nabla u^n_t \|_2^2 \,dt = 0$.
We need a strengthened version of this relation, which is presented in
the next corollary whose proof is easy, so you omit it.
\begin{corollary}\label{coefficientfn}
Let $f, f^n  \in\mathbf{L}^2 ([0,T] \times\mathbb{ R}^d ; \mathbb
{R } ), n \in\mathbb{ N},$ be such that\break
$ \lim_{n\to\infty} \int_0^T \| f^n_t - f_t \|_2^2 \,dt = 0$. Then
the solutions $(u^n)_{n\in\mathbb{N}}$
of the equations
\[
\bigl(\partial_t + \tfrac{1}{2} \Delta\bigr) u^{n} - n u^{n} + f^n = 0,
\]
with final condition $u_T^n = 0$, satisfy the relation $ \lim_{n\to
\infty} \int_0^T \|\nabla u^n_t \|_2^2 \,dt = 0$.
\end{corollary}
%
%
\begin{corollary}
\label{coefficientgn}
Let $g^n , g \in\mathbf{L}^2 ([0,T] \times\mathbb{ R}^d ; \mathbb
{R}^d )$ be such that\break
$ \lim_{n\to\infty} \int_0^T \| g^n_t - g_t \|_2^2 \,dt = 0$. Then
the solutions $(u^n)_{n\in\mathbb{N}}$
of the equations
\[
\bigl(\partial_t + \tfrac{1}{2} \Delta\bigr) u^{n} - n u^{n} + \operatorname{div} g^n =
0,
\]
with final condition $u_T^n = 0$, satisfy the relation $ \lim_{n\to
\infty} \int_0^T \|\nabla u^n_t \|_2^2 \,dt = 0$.
\end{corollary}

\begin{pf}
We regularize $g$ by setting $ g_{i,t}^{\epsilon} =
P_{\epsilon} g_{i,t}$ for $ i =1, \ldots,d, $ $\epsilon> 0$,
$t \in[0,T]$. Then $ g_{i}^{\epsilon} \in H^{1}_0 (\mathbb{R}^d)$
and $f^{\epsilon} = \operatorname{div} g^{\epsilon}$ is in $\mathbf{L}^2
([0,T] \times\mathbb{ R}^d ; \mathbb{R} )$. Moreover, we have
$ \lim_{\epsilon\to0} \int_0^T \| g^{\epsilon}_t - g_t \|_2^2\, dt =
0$. Let ${u}^{\epsilon,n}$ be the solution of the equation
\[
\bigl(\partial_t + \tfrac{1}{2} \Delta\bigr) u^{\epsilon,n} - n {u}^{\epsilon
,n} + f^{\epsilon} = 0,
\]
with final condition ${u}_T^{\epsilon,n} = 0$. By Lemma 5 of \cite
{Denis2}, one has
\begin{eqnarray*}
\int_0^T \|\nabla u^n_t - \nabla u^{\epsilon,n}_t \|_2^2 \,dt &\leq& c
\int_0^T \|g^n_t- g^{\epsilon}_t\|_2^2 \,dt\\
& \leq&
c \int_0^T (\|g^n_t- g_t\|_2^2 + \|g^{\epsilon}_t - g_t \|_2^2 ) \,dt.
\end{eqnarray*}
On the other hand, Lemma \ref{coefficientf} implies, for $ \epsilon$ fixed,
$ \lim_{n\to\infty} \int_0^T \|\nabla u^{\epsilon,n}_t \|_2^2 \,dt =
0$. From these facts, one easily concludes the proof.
\end{pf}

\begin{lemma}\label{coefficienthn}
Let $h , h^n, n \in\mathbb{N}$, be $\mathbf{L}^2(\mathbb{R}^d;
\mathbb{R}^{d^{1}})$-valued predictable processes on $ [0,T]$ with
respect to
$ (\mathcal{F}_{t,T}^B )_{t\geq0}$ and such that
\[
\mathbb{E} \int_0^T \|h_t\|_2^2 \,dt < \infty,\qquad
\mathbb{E} \int_0^T \| h_t^n\|_2^2 \,dt < \infty
\]
and
\[
\lim_{n\to\infty} \mathbb{E} \int_0^T \| h^n_t - h_t \|_2^2 \,dt = 0.
\]
Let $(u^n)_{n\in\mathbb{N}}$ be the solutions
of the equations
\[
d u_t^{n} + \bigl[ \tfrac{1}{2} \Delta u_t^{n} - n u_t^{n}\bigr]\, dt + h_t^n \cdot
\overleftarrow{dB}_t = 0,
\]
with final condition $u_T^n = 0$, for each $n\in\mathbb{N}$. Then one has
\[
 \lim_{n\to\infty} \int_0^T \|\nabla u^n_t \|_2^2 \,dt = 0.
 \]
\end{lemma}

\begin{pf}
We regularize the process $h$ by setting $ \bar
{h}_{i,t}^{\epsilon} = P_{\epsilon} h_{i,t}$ for $ i =1, \ldots
,\break d_1, \epsilon> 0$,
$t \in[0,T]$. Then $ \bar{h}_{i,t}^{\epsilon} \in H^{1}_0 (\mathbb{R}^d)$
and $ \mathbb{E} \int_0^T \|\nabla\bar{h}_t^{\epsilon} \|_2^2 \,dt <
\infty$ and\break
$ \lim_{\epsilon\to0} \mathbb{E} \int_0^T \| \bar
{h}_{t}^{\epsilon} - h_t \|_2^2 \,dt = 0$.
Let ${u}^{\epsilon,n}$ be the solution of the equation
\[
d{u}_t^{\epsilon,n} + \tfrac{1}{2} \Delta{u}_t^{\epsilon,n} - n
{u}_t^{\epsilon,n} + \bar{h}_t^{\epsilon}\cdot\overleftarrow{dB}_t
= 0
\]
with final condition ${u}_T^{\epsilon,n} = 0$, for each $n \in\mathbb
{N}$. The relation (iii) of Proposition 6 in \cite{Denis2} written
with respect to the Hilbert space $H=H_0^1 (\mathbb R^d)$ takes the form
\[
\mathbb{E} \biggl[ \|\nabla u_t^{\epsilon,n} \|_2^2 + \int_t^T \biggl\| \frac
{1}{2} \Delta u_s^{\epsilon,n} \biggr\|^2 \,ds + n
\int_t^T \|\nabla u_s^{\epsilon,n} \|_2^2 \,ds \biggr] = \mathbb{E } \int
_t^T \|\nabla\bar{h}{}^{\epsilon}_s \|_2^2 \,ds .
\]
In particular, one has
\[
\int_t^T \|\nabla u^{\epsilon,n}_s \|^2 \,ds \leq\frac{1}{n} \int
_t^T \|\nabla\bar{h}{}^{\epsilon}_s \|_2^2 \,ds.
\]
Now we write the relation (iii) of Proposition 6 in \cite{Denis2} for
the solution $ u^n - u^{\epsilon,n}$
with respect to the Hilbert space $H= \mathbf{L}^2(\mathbb R^d)$,
\begin{eqnarray*}
&&\mathbb{E} \biggl[ \|u_0^n- u_0^{\epsilon,n} \|^2 + \int_0^T
\|\nabla u_s^n- \nabla u_s^{\epsilon,n} \|_2^2 \,ds + n \int_0^T
\|u_s^n- u_s^{\epsilon,n} \|_2^2 \,ds \biggr]\\
&&\qquad = \mathbb{E } \int_0^T \|\bar
{h}_s^{n}-\bar{h}_s^{\epsilon} \|_2^2 \,ds .
\end{eqnarray*}
In particular, one obtains
\[
\mathbb{E} \int_0^T
\|\nabla u_s^n- \nabla u_s^{\epsilon,n} \|_2^2 \,ds \leq\mathbb{E }
\int_0^T \|\bar{h}_s^{n}-\bar{h}_s^{\epsilon} \|_2^2 \,ds.
\]

From this and the preceding inequality, one deduces
\[
\limsup_{n\to\infty} \mathbb{E} \int_0^T
\|\nabla u_s^n \|_2^2 \,ds \leq\mathbb{E } \int_0^T \|\bar{h}_s-\bar
{h}_s^{\epsilon} \|_2^2 \,ds.
\]
Letting $\epsilon\to0$, one deduces the relation from the statement.
\end{pf}

\begin{lemma}
\label{obstacleS}
Let $v \dvtx [0,T] \times\mathbb{R}^d \rightarrow\mathbb{R}$ be a
function such that the process $ ( v_t (W_t) )_{ t \in[0,T] }$ admits
a version $ S = (S_t )_{t \in[0,T]}$ with continuous trajectories on
$[0,T]$ and such that the random variable $ S^* = \sup_{0\leq t \leq
T} S_t$ satisfies the condition $ \mathbb{E}^m [ S^* ]^2 < \infty$.
Let $u^n$ be the solution of the equation
\[
\bigl(\partial_t + \tfrac{1}{2} \Delta\bigr) {u}^{n} - n {u}^{n} + n v = 0,
\]
with the terminal condition $u_T^n = v_T$. Let $Y^n = (Y_t^n )_{ t \in
[0,T]}$ be a continuous version of the process
$ (u_t^n (W_t) )_{t \in[0,T]}$, for each $n \in\mathbb{N}$. Then the
following holds:
\[
\lim_{ n \to\infty} \mathbb{E}^m \Bigl[
\sup_{ 0 \leq t \leq T} | Y_t^n - S_t |^2 \Bigr]= 0.
\]
\end{lemma}

\begin{pf}
Let us set $ \overline{u}^n_t = e^{-nt} u_t^n$ and
observe that this function is a solution of the equation
\[
\bigl(\partial_t + \tfrac{1}{2} \Delta\bigr) \overline{{u}}^{n} + \overline{ v } = 0,
\]
with $ \overline{v}_t = e^{-nt} v_t$ and terminal condition $\overline{u}_T^n =
\overline{v}_T$. Writing the representation of Theorem \ref{FK}
with $g=h=0$ for $ \overline{u}^n (W_t)$, one obtains
\[
\label{Ito:n}
\overline{u}_t^n ( t,W_{t} ) = e^{-nT} v_T - \sum_{i =1 }^d \int
_{t}^{T}\partial_{i}\overline{u}_r^n (W_{r} )
\,dW_{r}^{i}+ n \int_t^T e^{-nr} v_r (W_r)\, dr,
\]
and this leads to the representation of our process $Y^n,$ given by
\[
Y_t^n = \mathbb{E}^m \biggl[ e^{-n(T-t)} S_T + n \int_t^T e^{-n (r -t)}
S_r\,dr \Big| \mathcal{F}_t \biggr].
\]
Then one has
\[
|S_t - Y_t | \leq\mathbb{E}^m \biggl[ \biggl| S_t - e^{-n(T-t)} S_T - n \int
_t^T e^{-n (r -t)} S_r\,dr \biggr|
\Big| \mathcal{F}_t \biggr].
\]
Let us denote by
\[
V^n = \sup_{0\leq t \leq T} \biggl| S_t - e^{-n(T-t)} S_T - n \int_t^T
e^{-n (r -t)} S_r\,dr \biggr|.
\]
Obviously, one has $ V^n \leq2 S^*$. On the other hand, one has for
any fixed $\delta>0$,
\begin{equation}
\label{Vn}
V^n \leq\sup_{|t-s| \leq\delta} |S_t - S_s| + 2 e^{-n\delta} S^*.
\end{equation}
This follows from Lemma \ref{mc}. From the inequality \eqref{Vn}, one
deduces that $ \lim_{n\to\infty} V^n = 0, \mathbb{P}^m$-a.s., and hence from the dominated convergence theorem, one gets
$ \lim_{n\to\infty} \mathbb{E}^m [V^n ]^2= 0.$ Since
\[
|S_t - Y_t^n | \leq\mathbb{E}^m [V^n | \mathcal{F}_t ],
\]
Doob's theorem implies the assertion of the lemma.
\end{pf}

Finally, we mention the following calculus lemma.
\begin{lemma}
\label{mc}
Let $\varphi\in C ([0,1]; \mathbb{R} )$ and $\delta\in(0,T)$,
$\lambda> 0$. Then one has
\[
\biggl| \lambda\int_0^{\delta} e^{-\lambda t}\varphi(t) \,dt + e^{-\lambda
\delta}\varphi(\delta) - \varphi(0)\biggr|
\leq\sup_{0 \leq t \leq\delta} |\varphi(t) - \varphi(0)|
\]
and
\begin{eqnarray*}
&&\biggl| \lambda\int_t^{T} e^{-\lambda(s-t)}\varphi(s) \,ds + e^{-\lambda
(T-t)}\varphi(T) - \varphi(t) \biggr|\\
&&\qquad\leq\sup_{|s-r| \leq\delta, s\geq0} |\varphi(s) - \varphi(r)| +
2 e^{-\lambda\delta}\|\varphi\|_{\infty} .
\end{eqnarray*}
\end{lemma}

\begin{pf}
The first inequality follows from the relation
$ \lambda\int_0^{\delta}e^{- \lambda t} \,dt + e^{-\lambda\delta}
=1$. In order to check the second relation, one dominates the
expression of the left-hand side by
\begin{eqnarray*}
&&\biggl| \lambda\int_t^{t+\delta} e^{-\lambda(s-t)} \varphi(s) \,ds  +
e^{-\lambda\delta}\varphi(t+\delta) - \varphi(t) \biggr|\\
&&\qquad{} + e^{- \lambda\delta} \biggl|\lambda\int_{t+\delta}^{T} e^{-\lambda
(s-(t+\delta))} \varphi(s) \,ds + e^{-\lambda(T- (t+\delta))}\varphi
(T) - \varphi(t+\delta) \biggr|
\end{eqnarray*}
and then apply the first relation to dominate the first term.
\end{pf}

\begin{appendix}
\section*{Appendix}\label{app}
The next lemma is a classical result in convex analysis, known as
Mazur's theorem (see \cite{Brezis}, Remark 5, page 38).
We state here the result with some notation that is useful for our
proof. Let $X$ be a Banach space and $(x_n)_{n\in\mathbb{N}}$ a
sequence of elements in $X$. We call finite family of coefficients of a
convex combination a family $ a = \{ \alpha_i | i \in I\}$ where $I$
is a finite subset of $ \mathbb{N}$,
$\alpha_i > 0$ for each $ i \in I$ and $ \sum_{ i \in I} \alpha_i =
1$. The convex combination that corresponds to such a family of
coefficients is the point expressed in terms of our sequence by $ \sum
_{i \in I} \alpha_i x_i$.

\begin{lemma} \label{Mazur}
Let $(x_n)_{n\in\mathbb{N}}$ be a weakly convergent sequence of
elements in $X$ with limit $x$. Then there exits a sequence
$(a^k)_{k\in\mathbb{ N}}$ of families of coefficients of convex
combinations, $a^k =\{ \alpha_i^k | i \in I_k \}$, such
that the corresponding convex combinations $x^k = \sum_{i\in I_k}
\alpha_i^k x_i$, $k \in\mathbb{N},$ converge strongly to $x\dvt \lim_{ k \to\infty} \|x^k - x\| =0.$
\end{lemma}
\end{appendix}

%

\printaddresses

\end{document}